\documentclass[sn-mathphys-num]{sn-jnl}

\usepackage{graphicx}%
\usepackage{multirow}%
\usepackage{amsmath,amssymb,amsfonts}%
\usepackage{amsthm}%
\usepackage{mathrsfs}%
\usepackage[title]{appendix}%
\usepackage{xcolor}%
\usepackage{textcomp}%
\usepackage{manyfoot}%
\usepackage{booktabs}%
\usepackage{listings}%
\usepackage{bm}
\usepackage{subdepth}
\usepackage{enumitem}

\newtheorem{theorem}{Theorem}[section]
\newtheorem{proposition}[theorem]{Proposition}%
\newtheorem{lemma}[theorem]{Lemma}
\newtheorem{remark}[theorem]{Remark}%

\numberwithin{equation}{section}

\def\mcalB {{\mathcal{B}}}
\def\mcalC {{\mathcal{C}}}
\def\mcalD {{\mathcal{D}}}

\def\mcalP {{\mathcal{P}}}

\def\mcalS {{\mathcal{S}}}
\def\mcalX {{\mathcal{X}}}
\def\E     {{\mathbb{E}}}
\def\bR    {{\mathbb{R}}}

\def\dd    {{{\hskip 1pt}\rm{d}}}

\def\tr    {\mathop{\rm{tr}}}

\def\P    {\mathbb{P}}

\def\cl   {\mbox{\rm{cl}}}

\def\XX   {\bm{X}}
\def\TT   {\mb{T}}
\def\xx   {\bm{x}}
\def\ww   {\bm{w}}

\def\E    {\mathbb{E}}

\def \tt {\boldsymbol{t}}

\def \vv {\boldsymbol{v}}
\def \ww {\boldsymbol{w}}

\def \bzero {\boldsymbol0}

\def\balpha {{\boldsymbol{\alpha}}}

\def\bmu {{\boldsymbol{\mu}}}
\def\bnu {{\boldsymbol{\nu}}}
\def\bPsi {\boldsymbol{\Psi}}

\def\bLambda {{\boldsymbol{\Lambda}}}

\def\bSigma {{\boldsymbol{\Sigma}}}
\def\bXi    {{\boldsymbol{\Xi}}}
\def\bomega {{\boldsymbol{\omega}}}

\def\btheta {{\boldsymbol{\theta}}}
\def\bTheta {{\boldsymbol{\Theta}}}

\def \AA {\boldsymbol{A}}
\def \BB {\boldsymbol{B}}

\def \HH {\boldsymbol{H}}
\def \II {\boldsymbol{I}}

\def \MM {\boldsymbol{M}}

\def \SS {\boldsymbol{S}}
\def \TT {\boldsymbol{T}}
\def \UU {\boldsymbol{U}}

\def \WW {{\boldsymbol{W}}}
\def \XX {\boldsymbol{X}}

\def\e{{\mathbb E}}

\def\Cov{{\mathrm{Cov}}}

\DeclareMathOperator{\diag}{diag}

\begin{document}

\title[The Truncated Multivariate Normal Distributions]{Non-Steepness and Maximum Likelihood Estimation Properties of the Truncated Multivariate Normal Distributions}


\author*[1]{\fnm{Michael} \sur{Levine}}\email{mlevins@purdue.edu}

\author[2]{\fnm{Donald} \sur{Richards},}\email{dsr11@psu.edu}

\author[1]{\fnm{Jianxi} \sur{Su}}\email{jianxi@purdue.edu}

\affil*[1]{\orgdiv{Department of Statistics}, \orgname{Purdue University}, \orgaddress{\city{West Lafayette}, \postcode{IN 47906}, \country{U.S.A}}}

\affil[2]{\orgdiv{Department of Statistics}, \orgname{Pennsylvania State University}, \orgaddress{\city{University Park}, \postcode{PA 16802}, \country{U.S.A}}}


\abstract{This article considers exponential families of truncated multivariate normal distributions with one-sided truncation for some or all coordinates.  We observe that if all components are one-sided truncated then this family is not full.  The family of truncated multivariate normal distributions is extended to a full family, and the extended family is investigated in detail.  We identify the canonical parameter space of the extended family and  establish that the family is not regular and not even steep. We also consider maximum likelihood estimation for the location vector parameter and the positive definite (symmetric) matrix dispersion parameter of a truncated non-singular multivariate normal distribution.  It is shown that if the sample size is sufficiently large then, almost surely, the maximizer of
the likelihood function is unique, provided that it exists.  It is also shown that each solution to the score equations for the location and dispersion parameters satisfies the method-of-moments equations.  Finally, it is observed that similar results arise in the case of an arbitrary number of truncated components.}




\maketitle


\section{Introduction}
\label{sec_introduction}

In this article we continue our investigations in \cite{Levine} of the one-sided truncated multivariate normal distributions, a family of distributions which has been applied in numerous areas, including simultaneous equations inference and multivariate regression \cite{Amemiya}, economics and econometrics \cite{Heckman, Hong}, educational studies \cite{Jin,Levine}, agriculture \cite{Yu}, and inference with censored biomedical data \cite{Lachos,Lin}.  Here we are motivated, particularly by the results of Cohen \cite{Cohen}, del Castillo \cite{delCastillo}, and Hegde and Dahiya \cite{Hegde}, to study the problem of maximum likelihood estimation for these distributions.

It is well-known that many exponential families appearing in practice are regular and hence are steep.  For instance, the two-sided truncated univariate normal distributions, with support restricted to a finite interval, form a regular, and hence steep, exponential family (Barndorff-Nielsen \cite[p.~161]{Barndorff}).  There are also several exponential families that are steep but non-regular, an example being the family of inverse Gaussian distributions; see Bar-Lev and Enis \cite[p.~1516]{Bar_Lev}, Barndorff-Nielsen \cite[p.~117]{Barndorff}, and Brown \cite[p.~72]{Brown}.  Examples of non-steep multivariate exponential families arising in applied research are provided by Sundberg \cite[p.~251 and p.~259]{Sundberg}, who studied the Strauss model and inference for data on refindings of ringed birds; Bedbur, Kamps, and Imm \cite{Bedbur_etal}, on maximum likelihood inference for the parameters of the Conway-Maxwell-Poisson distribution; and Bar-Lev and Ridder \cite{Bar_Lev_Ridder}, for non-steep exponential models for count data.  

In the univariate case, del Castillo \cite{delCastillo} extended the one-sided truncated normal distributions to a full exponential family of distributions, and proved that the extended family is not steep and therefore is not regular.  That the extended family of distributions is not steep implies that the likelihood equations may have no solutions, although the likelihood function may have its maximum on the boundary of the parameter space \cite[Corollary 4.4]{delCastillo}.  Therefore del Castillo's results have deeper and profound implications for statistical inference for the parameters of the one-sided truncated normal distributions.

During the three decades that have elapsed since the appearance of del Castillo's article \cite{delCastillo}, there has been no subsequent investigation of the multivariate truncated normal distributions as an exponential family.  It is this observation of the absence of such results that initially motivated our research here and in \cite{Levine}.  We consider in this article the family, $\mcalP_\mcalD$, of extended truncated multivariate $d$-dimensional normal distributions such that every univariate component is one-sided truncated.  Our first main result, given in Section \ref{sec_non_regularity} is that the family $\mcalP_\mcalD$ is a multivariate exponential family.  We prove that the canonical (or natural) parameter space, $\mcalD$, consists of a union of $d+1$ subsets, all but one of which are not open; consequently, we deduce that $\mcalD$ is not open and therefore $\mcalP_\mcalD$ is not regular.

In Section \ref{sec_non_steep}, we prove that the family $\mcalP_\mcalD$ is non-steep.  As the marginal distributions of the truncated multivariate normal distributions generally are not truncated normal \cite{Levine}, the non-steepness property in the multivariate case does not follow from the univariate case.  By analyzing the cumulant-generating function of the family $\mcalP_\mcalD$, we prove that the gradient of that function remains bounded near certain boundary points in $\mcalD$, hence we find that $\mcalP_\mcalD$ is not steep.

We consider in Section \ref{sec_mle} the problem of maximum likelihood estimation for $\bmu$, the location vector parameter, and $\bSigma$, the positive definite (symmetric) matrix dispersion parameter, of a truncated non-singular multivariate normal distribution.  We prove that the solutions to the score equations also satisfy the method-of-moments equations, thereby generalizing results of various authors in the univariate, bivariate, and trivariate cases (see Cohen \cite[Sections 2.3 and 12.2]{Cohen}).  As a consequence of these results we obtain a necessary condition for the existence of solutions to the score equations for ($\bmu,\bSigma)$, generalizing a result of del Castillo \cite[Theorem 4.2]{delCastillo}. Section \ref{sec_several_untruncated} considers the case where only one coordinate has been truncated. We show that the resulting exponential family is also not steep. We also note that this result can be extended to the case where more than one, but not all, coordinates are truncated. Section \ref{sec_conclusions} summarizes the results presented in this manuscript and suggests some avenues for future research. Finally, Section \ref{sec_proofs} contains proofs of the results presented here.

\section{Non-regularity properties of the family \texorpdfstring{$\mathcal{\mcalP_\mcalD}$}{Sigma}}\label{sec_non_regularity}

In the sequel, we write $\xx > \bzero$ to mean that $\xx \in \bR^d_+$ and write $\xx < \bzero$ whenever $-\xx > \bzero$; more generally, we write $\xx > \xx_0$ if $\xx - \xx_0 > \bzero$.

The exponential family considered in the present article arises from a random vector $\XX = (X_1,\ldots,X_d)' \in \bR^d$, having a truncated multivariate normal distribution such that each component of $\XX$ is singly-truncated, i.e., truncated below or above, at a specified finite point.  By applying sign changes to the components, if necessary, we may assume that each $X_j$ is truncated below at a given value $a_j$, $j=1,\ldots,d$.  Also, by the change-of-variables $x_j \to x_j + a_j$, $j=1,\ldots,d$, we may assume, without loss of generality, that $a_j = 0$ for all $j=1,\ldots,d$, so that the distribution of $\XX$ is supported on $\bR_+^d$, the \textit{nonnegative} orthant.  Here, we are using the property that the family of truncated multivariate normal distributions is invariant under changes of sign or location; see Levine, \textit{et al.} \cite{Levine}.  

For a location parameter $\bmu \in \bR^d$, and a positive definite (symmetric) $d \times d$ matrix dispersion parameter $\bSigma$, let $\XX = (X_1,\ldots,X_d)' \in \bR^d$ be a random vector having the probability density function
\begin{equation}
\label{eq_trunc_normal_pdf}
f(\xx;\bmu,\bSigma) = [C(\bmu,\bSigma)]^{-1} \, \exp\big[- \tfrac12(\xx-\bmu)'\bSigma^{-1}(\xx-\bmu)\big],
\end{equation}
$\xx > \bzero$, where
\begin{equation}
\label{eq_normaliz_const}
C(\bmu,\bSigma) = \int_{\ww > \bzero} \exp\big[-\tfrac12 (\ww-\bmu)'\bSigma^{-1}(\ww-\bmu)\big] \dd \ww
\end{equation}
is the corresponding normalizing constant.  The truncated multivariate normal distributions given by \eqref{eq_trunc_normal_pdf} are the subject of the present article.  

Let $\btheta = \bSigma^{-1} \bmu$ and $\bTheta = \tfrac12 \bSigma^{-1}$ denote the \textit{canonical parameters} of the family of truncated multivariate normal distributions.  Then the density function \eqref{eq_trunc_normal_pdf} is proportional to $\exp(\btheta'\xx - \xx'\bTheta\xx)$, $\xx > \bzero$, and this leads to the formulation of these distributions as an exponential family.  

The multivariate model \eqref{eq_trunc_normal_pdf} is a family of exponential type, with an open parameter space when expressed in terms of $(\bmu,\bSigma)$, but it is not a full exponential family.  We will show in due course that the full family includes some boundary models, which implies that $\mcalD$, the corresponding canonical parameter space, is not open.  

We now provide a brief review of some classical concepts in the theory of exponential families, concepts that are treated in detail by Barndorff-Nielsen \cite{Barndorff}, Brown \cite{Brown}, Sundberg \cite{Sundberg}, and other authors.  For a positive integer $k$, let $\mcalX \subseteq \bR^k$ be a measurable space.  Also let $\beta$ be a sigma-finite, positive measure on $\mcalX$ such that $\beta$ is not concentrated on an affine subspace.  Let $T : \mcalX \to \bR^k$ be a statistic, and denote by $\beta_T$ the marginal distribution of $T$.  For a $k \times 1$ column vector $\tt = (t_1,\ldots,t_k)' \in \bR^k$ and a parameter $\bomega = (\omega_1,\ldots,\omega_k)' \in \bR^k$, define the Laplace transform,
\begin{equation}
\label{eq_LT_omega}
L_{T}(\bomega) = \int_{\mcalX} \exp(\bomega'T) \dd\beta(T) = \int_{\bR^k} \exp(\bomega'\tt) \dd\beta_{T}(\tt).
\end{equation}
The set $\mcalD = \{\bomega \in \bR^k: L_{T}(\bomega) < \infty\}$ is called the \textit{canonical parameter space}.  By applying H\"{o}lder's inequality, we find that $\mcalD$ is convex.  Further, define the \textit{cumulant-generating function}, $K_{T}(\bomega) = \log L_{T}(\bomega)$, $\bomega \in \mcalD$; then $K_T$ is strictly convex whenever it is non-degenerate.

For each $\bomega \in \mcalD$, define $p(\tt;\bomega) = \exp\big(\bomega'T(\tt) - K_{T}(\bomega)\big)$, $\tt \in \bR^k$; then the measure $p(\tt;\bomega) \dd\beta(\tt)$, $\tt \in \bR^k$, is a probability distribution.  The family of distributions $\mcalP_\mcalD = \{p(\tt;\bomega) \dd\beta(\tt), \tt \in \bR^k: \bomega \in \mcalD\}$ is an \textit{exponential family}, and the family is called \textit{full} if the parameter $\bomega$ is allowed to vary over the entire space $\mcalD$ \cite[p.~2]{Brown}.

The family $\mcalP_\mcalD$ is said to be \textit{regular} if $\mcalD$ is an open set.  If $\mcalD$ is not open then denote by $\Omega$ the set $\hbox{Int}(\mcalD)$, the interior of $\mcalD$.  We consider only the non-trivial case in which $\Omega$ is non-empty, and then we define ${\cal P}_{\Omega} = \{p(\tt;\bomega) \dd\beta(\tt), \tt \in \bR^k: \bomega \in \Omega\}$.

Starting with \eqref{eq_LT_omega}, let $\mcalX = \bR^d_+$ and let $\beta$ be the Lebesgue measure on $\bR_+^d$, written $\textrm{d}\beta(\tt) \equiv \textrm{d}\tt$.  Choose as sufficient statistics the collection
\begin{equation}
\label{eq_T_of_tt}
T(\tt) = \{t_i: i=1,\ldots,d\} \cup \{t_i t_j: 1 \le i \le j \le d\},
\end{equation}
so that $T$ is a statistic of dimension 
$$
k = d + \tfrac12 d(d+1).
$$
The Laplace transform, $L_T$, of the resulting distribution $\beta_T$ is the Laplace transform, with respect to $\textrm{d}\tt$, of the set of statistics $T(\tt)$.

Denote by $\mcalS^{d \times d}$ the vector space of $d \times d$ symmetric matrices and observe that each linear combination of the sufficient statistics in $T$ is of the form $\btheta'\tt - \tt' \bTheta \tt$ for some $\btheta \in \bR^d$ and $\bTheta \in \mcalS^{d \times d}$.  Therefore $L_T$ is parametrized by the pairs $\bomega = (\btheta,\bTheta) \in \bR^d \times \mcalS^{d \times d}$, the Laplace transform of $\beta_T$ is
\begin{equation}
\label{eq_Laplace_transform}
L_T(\btheta,\bTheta) = \int_{\ww > \bzero} \exp(\btheta'\ww - \ww' \bTheta \ww) \, \dd\ww,
\end{equation}
whenever the integral converges, and the corresponding cumulant-generating function is 
\begin{equation}
\label{eq_KT_cgf}
K_T(\btheta,\bTheta) = \log L_T(\btheta,\bTheta) = \log \int_{\ww > \bzero} \exp(\btheta'\ww - \ww' \bTheta \ww) \, \dd\ww.
\end{equation}

It is now evident that, subject to convergence of the integral \eqref{eq_Laplace_transform}, the 
function 
\begin{equation}
\label{eq_trunc_normal_nef_pdf}
p(\tt;\btheta,\bTheta) = \exp\big(\btheta'\tt - \tt' \bTheta \tt - K_T(\btheta,\bTheta)\big), \qquad \tt \in \bR_+^d,
\end{equation}
is a probability density function.  Note that the function \eqref{eq_trunc_normal_nef_pdf} may be a probability density function even if $\bTheta$ is not positive definite, e.g., if $\btheta < \bzero$ and $\bTheta = \bzero$.   Consequently, we view the class of densities \eqref{eq_trunc_normal_nef_pdf} as representing a family of \textit{extended} truncated multivariate normal distributions.  

To determine the resulting exponential family $\mcalP_\mcalD$, we identify the corresponding canonical parameter space,
\begin{equation}
\label{eq_D_space}
\mcalD = \{(\btheta,\bTheta) \in \bR^d \times \mcalS^{d \times d}: L_T(\btheta,\bTheta) < \infty\}.
\end{equation}
Denote by $\bzero$ any zero vector or matrix, with its dimension determined by the context.  For $0 \le r \le d$, denote by $\mcalS_r^{d \times d}$ the set of $d \times d$ positive semidefinite (symmetric) matrices of rank $r$; thus, for $r \le d-1$ the matrices in $\mcalS_r^{d \times d}$ are singular, and the matrices in $\mcalS_d^{d \times d}$ are positive definite and hence non-singular.  

In proving the following result, we encounter a collection of parameter spaces $\mcalD_r$, $0 \le r \le d$, that are defined in \eqref{eq_A_0}, \eqref{eq_A_r}, and \eqref{eq_A_d}.  We will have $\mcalD_0 = \{(\btheta,\bzero): \btheta < \bzero\}$ and $\mcalD_d = \{(\btheta,\bTheta): \theta \in \bR^d, \bTheta \in \mcalS_d^{d \times d}\}$, and the remaining ``boundary spaces'' $\mcalD_r$, $1 \le r \le d-1$, are defined by technical constructions during the course of the proof.  

\medskip

\begin{theorem}
\label{thm_parameter_space}
The family, $\mcalP_\mcalD$, of extended truncated multivariate normal distributions constitutes a full and minimal exponential family with canonical parameter space
\begin{equation}
\label{eq_D_nat_param_space}
\mcalD = \bigcup_{r=0}^d \mcalD_r.
\end{equation}
Further, the parameter spaces $\mcalD_0,\ldots,\mcalD_{d-1}$ are not open and $\mcalP_\mcalD$ is not regular.  
\end{theorem}

\bigskip

The proof of this result is given in subsection \ref{subsec_sec_non_regularity}.

\smallskip

For the case in which $d=1$, it follows by \eqref{eq_D_nat_param_space} that $\mcalD = \mcalD_0 \cup \mcalD_1$, where $\mcalD_0 = \{(\theta,0): \theta < 0\}$ and $\mcalD_1 = \{(\theta,\Theta): -\infty < \theta < \infty, \Theta > 0\}$.  This result agrees with the descriptions of $\mcalD$ given by del Castillo \cite{delCastillo} and Hegde and Dahiya \cite{Hegde}.

\section{Non-steepness properties of the family \texorpdfstring{$\mathcal{\mathbf{\mcalP_\mcalD}}$}{Sigma}}
\label{sec_non_steep}

In this section we show that the exponential family $\mcalP_\mcalD$ of extended truncated multivariate normal distributions is not steep.

Let $\mcalP$ be a generic natural exponential family on $\bR^k$ with canonical parameter space $\mcalD$ and cumulant-generating function $K$.  Let $\Omega = \hbox{Int}(\mcalD)$ be the interior of $\mcalD$, $\Omega^c$ be the complement of $\Omega$, and $\cl(\Omega)$ and $\cl(\Omega^c)$ be the closure of $\Omega$ and $\Omega^c$, respectively.  The \textit{boundary} of $\Omega$ is the set
\begin{equation}
\label{eq_Omega_boundary}
\partial\Omega = \cl(\Omega) \cap \cl(\Omega^c)
\end{equation}
\cite[p.~28, Exercise 8]{Gaal}, \cite[p.~336, Definition 5.5.3]{Horn}.

For $\bomega = (\omega_1,\ldots,\omega_k)\in \bR^k$, let $\nabla_{\bomega} = (\partial/\partial \omega_1,\ldots,\partial/\partial \omega_k)'$ be the gradient operator on $\bR^k$, and let $\|\cdot\|$ denote the Euclidean norm on $\bR^k$.  The family $\mcalP$ is said to be \textit{steep at the point} $\bomega_0 \in \partial\Omega$ if
\begin{equation}
\label{eq_steep}
\|\nabla_\bomega K(\bomega)\| \to \infty \ \ \textrm{as} \ \ \bomega \to \bomega_0.
\end{equation}
Further, $\mcalP$ is called \textit{steep} if it is steep at every point $\bomega_0 \in \partial\Omega$; see \cite[p.~86]{Barndorff}, \cite[p.~71]{Brown}.  It is well known that if $\mcalP$ is regular then it is steep \cite[p.~71]{Brown}.

We now apply the foregoing concepts to $\mcalP_\mcalD$, the family of truncated multivariate normal distributions.
By Theorem \ref{thm_parameter_space}, the canonical parameter space corresponding to $\mcalP_\mcalD$ is the set $\mcalD = \bigcup_{r=0}^d \mcalD_r$.  Given $\bomega_1 = (\btheta_1,\bTheta_1)$ and $\bomega_2 = (\btheta_2,\bTheta_2)$ in $\mcalD$, define the inner product,
$$
\langle \bomega_1,\bomega_2\rangle = \btheta_1'\btheta_2 + \tr(\bTheta_1\bTheta_2).
$$
This inner product induces for each $\bomega = (\btheta,\bTheta) \in \mcalD$ the norm,
\begin{equation}
\label{eq_D_norm}
\|\bomega\| = \langle \bomega,\bomega\rangle^{1/2} = \big[\|\btheta\|^2 + \|\bTheta\|_F^2)\big]^{1/2},
\end{equation}
where $\|\bTheta\|_F := [\tr(\bTheta^2)]^{1/2}$ is the \textit{Frobenius norm} of $\bTheta$ \cite[p.~291]{Horn}.

As defined earlier, for $\btheta = (\theta_1,\ldots,\theta_d)' \in \bR^d$, the gradient operator is $\nabla_\btheta = (\partial/\partial\theta_1,\ldots,\partial/\partial\theta_d)'$.
Also use the Kronecker symbol notation, $\delta_{i,j} = 1$ or $0$ according as $i=j$ or $i \neq j$, respectively; then for $\bTheta = (\theta_{i,j}) \in \mcalS^{d \times d}$, the natural gradient operator on $\mcalS^{d \times d}$ is the $d \times d$ symmetric matrix of differential operators 
$\nabla_\bTheta = \big(\tfrac12(1+\delta_{i,j})\partial/\partial \theta_{i,j}\big)$.

For $\bomega = (\btheta,\bTheta) \in \hbox{Int}(\mcalD)$, define the gradient operator $\nabla_\bomega = (\nabla_\btheta,\nabla_\bTheta)$.  Let $K:\bR^d \times \mcalS^{d \times d} \to \bR$ be such that all partial derivatives of $K$ exist and are continuous.  Consistent with \eqref{eq_D_norm}, we define
\begin{equation}
\label{eq_nabla_omega}
\|\nabla_\bomega K(\btheta,\bTheta)\| := \big[\|\nabla_\btheta K(\btheta,\bTheta)\|^2 + \|\nabla_\bTheta K(\btheta,\bTheta)\|_F^2\big]^{1/2}.
\end{equation}

We shall show that there exists $(\btheta_0,\bTheta_0) \in \partial\Omega$, and a sequence $(\btheta_n,\bTheta_n) \in \hbox{Int}(\mcalD)$, $n =1,2,3,\ldots$, such that $(\btheta_n,\bTheta_n) \to (\btheta_0,\bTheta_0)$ and $\|\nabla_\bomega K(\btheta_n,\bTheta_n)\|^2 \not\to \infty$ as $n \to \infty$.  Consequently, we will obtain the following result, the proof of which is given in subsection \ref{subsec_sec_non_steep}.

\medskip

\begin{theorem}
\label{thm_non_steep}
The exponential family $\mcalP_\mcalD$ is not steep.
\end{theorem}

\medskip

As a consequence of this result, it also follows that the family $\mcalP_\mcalD$ is not regular; cf., Brown \cite[p.~71]{Brown}.

\section{Maximum likelihood estimation}
\label{sec_mle}

Let $\xx_1,\ldots,\xx_n$ be a random sample from $\XX$, a random vector having the $d$-dimensional truncated normal distribution, with the probability density function given in \eqref{eq_trunc_normal_pdf} and the normalizing constant \eqref{eq_normaliz_const}.  As shown in Section \ref{sec_non_steep}, the family of extended truncated multivariate normal distributions is not steep and hence is not regular; therefore $(\widehat \bmu,\widehat \bSigma)$, the maximum likelihood estimator of $(\bmu,\bSigma)$ may not be a solution to the likelihood equations. In this section, we show that if $n > d$ then, with probability one, there exists on $\mcalD$ a maximum for the likelihood function.  We shall also consider the case in which $\bSigma$ is positive definite, and in that case we derive the score equations and obtain a necessary condition for the existence of solutions to the likelihood equations.

Let ``$\succcurlyeq$'' denote the L{\"o}wner (positive semidefinite) ordering on $\mcalS^{d \times d}$, the space of $d \times d$ symmetric matrices, i.e., $\AA \succcurlyeq \BB$ (or $\BB \preccurlyeq \AA$) if $\AA - \BB$  is positive semidefinite.  
Given a random sample $\xx_1,\ldots,\xx_n$, construct the statistics $\TT_1 = \overline{\xx} = n^{-1} \sum_{j=1}^n \xx_j$ and $\TT_2 = n^{-1} \sum_{j=1}^n \xx_j \xx_j'$.  Further, define the parameter set 
$$
\mcalC_d = \{(\btheta,\bTheta) \in \bR_+^d \times \mcalS^{d \times d}: \bTheta \succcurlyeq \btheta\btheta'\}.
$$
It is evident that $\mcalC_d$ is non-empty, 
and we show later that $\mcalC_d$ also is closed, convex, and contains the \textit{convex support}, i.e., the closed convex hull of the support, of $\beta_{(\TT_1,\TT_2)}$.

In the next theorem, which generalizes a result of del Castillo \cite[Section 4]{delCastillo}, it will also be shown that $(\TT_1,\TT_2)$ is sufficient for $(\bmu,\bSigma)$, and that if $n > d$ then $(\widehat\bmu,\widehat\bSigma)$ exists almost surely and is unique.  

\medskip

\begin{theorem}
\label{thm_existence_mles}
Let $\xx_1,\ldots,\xx_n$ be a random sample from $\XX$.  Then
\begin{enumerate}[label=(\roman*), align=parleft] 
\item
The pair of statistics $(\TT_1,\TT_2)$ is sufficient for $(\bmu,\bSigma)$.

\item
The set $\mcalC_d$ is closed, convex, and contains the convex support of $\beta_{(\TT_1,\TT_2)}$.

\item
If $n > d$ then, with probability one, $(\TT_1,\TT_2) \in {\rm{Int}}(\mcalC_d)$ and if the likelihood function has a maximizer on $\mcalD$ then that maximizer is unique.
\end{enumerate}
\end{theorem}

\medskip

\begin{remark}
\label{rem_existence_mles}
{\rm 
(i) In Theorem \ref{thm_existence_mles}(ii), it remains unknown whether the convex support of $\beta_{(\TT_1,\TT_2)}$ equals $\mcalC_d$.  We note that, since $\TT_1$ is the sample mean vector and $\TT_2 - \TT_1\TT_1'$ is a constant multiple of the sample covariance matrix then, in the untruncated case, $\TT_1$ and $\TT_2 - \TT_1\TT_1'$ are mutually independent and the support of $(\TT_1,\TT_2)$ is obtained using the explicit formulas for the marginal distributions of $\TT_1$ and $\TT_2 - \TT_1\TT_1'$.  In the truncated case, however, $\TT_1$ and $\TT_2 - \TT_1\TT_1'$ are not independent and there is no explicit description of their joint distribution, so it is accordingly more difficult to characterize the support.  

(ii) Concerning Theorem \ref{thm_existence_mles}(iii), the proof will follow from the fact that the family $\mcalP_\mcalD$ is minimal.  Since $\mcalP_\mcalD$ also is non-steep then it is possible that the maximum likelihood estimator of $(\bmu,\bSigma)$ may lie on the boundary of the natural parameter space; however, if the maximum likelihood estimator exists then it is uniquely determined since the family is minimal and full, and therefore the log-likelihood function is strictly concave on the convex natural parameter space.  In this regard, we also refer to the extensive discussion of Bedbur, \textit{et al.}~\cite[Section 3]{Bedbur_etal} who encountered a full, non-steep, exponential family for which maximum likelihood estimators can occur on the boundary.  
}\end{remark}

\medskip

We write $\XX \sim N_d(\bmu,\bSigma;\bzero)$ whenever $\XX$ has the density function \eqref{eq_trunc_normal_pdf}.  Also, we now assume that $\bSigma$ is positive definite.

Unlike the untruncated normal distribution, $\bmu$ and $\bSigma$ are not the mean and covariance matrix, respectively, of $\XX$, so we define
$$
\bnu = (\nu_1,\ldots,\nu_d) := \E(\XX)
$$
and
$$
\bLambda = (\lambda_{i,j}) := \Cov(\XX) \equiv \E[(\XX-\bnu)(\XX-\bnu)'].
$$
We note that research on the moments of the truncated multivariate normal distributions was initiated in \cite{Rosenbaum} and continued later in \cite{Amemiya,Leppard,Tallis}.  Among later publications, \cite{Kan} derived a moment-generating function for the truncated multivariate normal distribution and obtained recurrence relations for related integrals that involve the density function of the untruncated multivariate normal distributions; those relations resulted in useful explicit expressions for low-order moments of folded and truncated multivariate normal distributions.

Let us introduce the notation
$$
(\nabla C)(\bmu,\bSigma) := \nabla_{\tt}\, C(\tt,\bSigma)\Big|_{\tt = \bmu}
$$
for the gradient of $C(\bmu,\bSigma)$ with respect to $\bmu$ and
$$
(\nabla \nabla' C)(\bmu,\bSigma) := \nabla_{\tt}\phantom{'} \nabla_{\tt}' C(\tt,\bSigma)\Big|_{\tt = \bmu}
$$
for the corresponding Hessian matrix.  The following auxiliary result is useful as it leads to expressions for the mean and covariance matrix of the multivariate truncated normal distribution in terms of the corresponding parameters of the underlying multivariate normal distribution. These results enable us, in turn, to obtain convenient representations of the maximum likelihood estimators of parameters of the underlying multivariate normal distributions in terms of observations generated by a multivariate truncated normal distribution, provided that these estimators exist and are in the interior of the parameter space. 

\medskip

\begin{lemma}
\label{lem_mean_covariance}
Suppose that $\XX \sim N_d(\bmu,\bSigma;\bzero)$.  Then
\begin{equation}
\label{eq_nu_ML}
\bnu = \bmu + [C(\bmu,\bSigma)]^{-1} \ \bSigma \, (\nabla C)(\bmu,\bSigma)
\end{equation}
and
\begin{equation}
\label{eq_Lambda_ML}
\bLambda = \bSigma + [C(\bmu,\bSigma)]^{-1} \bSigma (\nabla \nabla' C)(\bmu,\bSigma) \bSigma - (\bnu - \bmu)(\bnu - \bmu)'.
\end{equation}
\end{lemma}

\medskip

Suppose that $(\widehat \bmu,\widehat \bSigma)$ exists.  By \eqref{eq_nu_ML}, \eqref{eq_Lambda_ML}, and the invariance principle of maximum likelihood estimation \cite[pp.~70--71]{Anderson}, we find that the pair $(\widehat \bmu,\widehat \bSigma)$ is related to $(\widehat\bnu,\widehat\bLambda)$, the maximum likelihood estimator of $(\bnu,\bLambda)$, through the equations
\begin{equation}
\label{eq_nu_MLE}
\widehat\bnu = \widehat\bmu + [C(\widehat\bmu,\widehat\bSigma)]^{-1} \ \widehat\bSigma \, (\nabla C)(\widehat\bmu,\widehat\bSigma)
\end{equation}
and
\begin{equation}
\label{eq_Lambda_MLE}
\widehat\bLambda = \widehat\bSigma + [C(\widehat\bmu,\widehat\bSigma)]^{-1} \widehat\bSigma \, (\nabla \nabla' C)(\widehat\bmu,\widehat\bSigma) \, \widehat\bSigma - (\widehat\bnu - \widehat\bmu)(\widehat\bnu - \widehat\bmu)'.
\end{equation}
Our next goal is to obtain the score equations for $(\widehat\bnu,\widehat\bLambda)$.   

In the sequel, for any vector $\balpha \in \mathbb{R}^d$, define
\begin{equation}
\label{eq_def_A}
\SS(\balpha) = \frac{1}{n} \sum_{j=1}^n (\xx_j-\balpha)(\xx_j-\balpha)'.
\end{equation}
In particular,
$$
\SS(\overline{\xx}) = \frac{1}{n} \sum_{j=1}^n (\xx_j - \overline{\xx})(\xx_j - \overline{\xx})'
$$
is the maximum likelihood estimator of the covariance  matrix in the classical untruncated setting.

\medskip

\begin{theorem}
\label{thm_score_equations}
Let $\xx_1,\ldots,\xx_n$ be a random sample from $\XX$ and suppose that $(\widehat \bmu,\widehat \bSigma)$, the maximum likelihood estimator of $(\bmu,\bSigma)$, exists and is in the interior of the parameter space.  Then $(\widehat \bmu,\widehat \bSigma)$ satisfies the equations
\begin{equation}
\label{eq_score_eqs_nu}
\widehat\bmu + [C(\widehat\bmu,\widehat\bSigma)]^{-1} \ \widehat\bSigma \, (\nabla C)(\widehat\bmu,\widehat\bSigma) = \overline{\xx},
\end{equation}
and
\begin{equation}
\label{eq_score_eqs_Lambda}
\widehat\bSigma + [C(\widehat\bmu,\widehat\bSigma)]^{-1} \widehat\bSigma \, (\nabla \nabla' C)(\widehat\bmu,\widehat\bSigma) \, \widehat\bSigma - (\overline{\xx} - \widehat\bmu)(\overline{\xx} - \widehat\bmu)' = \SS(\overline{\xx}).
\end{equation}
Moreover, $\widehat\bnu = \overline{\xx}$ and $\widehat\bLambda = \SS(\overline{\xx})$.  
\end{theorem}

\medskip

\begin{remark}
\label{rem_last}
{\rm
(i) In the untruncated case, the normalizing constant $C(\bmu,\bSigma)$ does not depend on $\bmu$ and therefore $\nabla_\bmu C(\bmu,\bSigma) \equiv \bzero$.  Then \eqref{eq_score_eqs_nu} and \eqref{eq_score_eqs_Lambda} reduce to $\widehat{\bnu} = \widehat{\bmu} = \overline{\xx}$ and $\widehat{\hbox{Cov}}(\XX) = n^{-1} \sum_{j=1}^n (\xx_j-\overline{\xx})(\xx_j-\overline{\xx})'$, which are well-known results for the untruncated normal distributions.

(ii) It also follows from \eqref{eq_score_eqs_nu} and \eqref{eq_score_eqs_Lambda} that the method-of-moments estimator of $(\bnu,\bLambda)$ for the truncated multivariate normal distribution is the same as in the untruncated case.  This result was noted by Cohen \cite[Sections 2.3 and 12.2]{Cohen} for the univariate truncated normal distributions and for multivariate normal distributions with a \textit{single truncated component variable}.  Therefore, Theorem \ref{thm_score_equations} extends Cohen's observation to the general multivariate setting in which all components are truncated.

We also remark that Cohen's result for the univariate truncated normal distributions or the multivariate normal distributions with one truncated component cannot be extended to the general truncated multivariate setting by using conditional distributions, for in the case of the truncated normal distributions, the conditional distributions generally are not truncated normal.

(iii) For the case in which only one component is truncated, Cohen derived explicit formulas for the maximum likelihood estimator.  However, it appears to be infeasible to obtain similar formulas if two or more components are truncated; moreover, Cohen did not consider such cases.
}\end{remark}

\smallskip

Due to the implicit nature of the score equations \eqref{eq_score_eqs_nu} and \eqref{eq_score_eqs_Lambda}, it is a challenging problem to derive necessary and sufficient conditions for the existence and uniqueness of solutions to those equations.  The following result provides a necessary condition for the existence of solutions to the score equations, and we also note that this condition was proved by del Castillo \cite[Theorem 4.2]{delCastillo} to be necessary and sufficient in the one-dimensional case.

\medskip

\begin{proposition}
\label{prop_necessary_condition}
Suppose $n > d$ and that there exists a solution $(\widehat\bmu,\widehat\SS)$ to the score equations \eqref{eq_score_eqs_nu} and \eqref{eq_score_eqs_Lambda}.  Then 
\begin{equation}
\label{eq_necessary}
0 \le (\overline{\xx} - \widehat\bmu)' \widehat\SS^{-1} (\overline{\xx} - \widehat\bmu) < 1.
\end{equation}
\end{proposition}

\smallskip

In the one-dimensional case in which $X$ has a normal distribution that is truncated to the left at $0$, it is shown in \cite{delCastillo,Hegde} that the domain of means, $\{(\E(X),\E(X^2))\}$, \textit{equals} the region $\{(x,y) \in \bR_+^2: x^2 \le y \le 2x^2\}$; equivalently, 
\begin{equation}
\label{eq_domain_means}
[\E(X)]^2 \le \E(X^2) \le 2 [\E(X)]^2.
\end{equation}
Moreover, the sufficiency of \eqref{eq_necessary} in the one-dimensional case stems from an upper bound on the domain of means.  

In the multivariate case, it was shown in Theorem \ref{thm_parameter_space} that $\E(\XX)\E(\XX') \preccurlyeq \E(\XX \XX')$, and this reduces in the one-dimensional case to the lower bound in \eqref{eq_domain_means}.  So, it remains to obtain an upper bound for $\E(\XX \XX') - \E(\XX)\E(\XX')$ in the multivariate case, and we do so in the following result.  Using the notation $\II_d$ for the identity matrix of order $d$, we have

\smallskip

\begin{proposition}
\label{prop_domain_means_upper}
Suppose that $\XX$ has the truncated multivariate normal distribution \eqref{eq_trunc_normal_pdf} with mean $\bnu$ and covariance matrix $\bLambda$.  Then 
\begin{equation}
\label{eq_domain_means_upper_1}
\E(\XX\XX') \preccurlyeq \E(\XX)\E(\XX)' + (\E (\|\XX\|^2) + \|\E(\XX)\|^2) \II_d,
\end{equation}
and the domain of means $\{(\E(\XX),\E(\XX\XX'))\}$ is contained in the set 
\begin{equation}
\label{eq_domain_means_upper_2}
\{(\btheta,\bTheta): \btheta\btheta' \preccurlyeq \bTheta \preccurlyeq \btheta\btheta' + (\tr\bTheta + \|\btheta\|^2) \II_d\}.
\end{equation}
\end{proposition}

\smallskip

We remark that although \eqref{eq_domain_means_upper_2} is vacuous in the case $d=1$, it suggests that any necessary and sufficient conditions for the existence of a solution $(\widehat\bmu,\widehat\SS)$ to the score equations \eqref{eq_score_eqs_nu} and \eqref{eq_score_eqs_Lambda} will involve not only $\bar{\xx}$, $\widehat\SS$, but also polynomials in $\|\bar{\xx}\|$ and $\widehat\SS$.  

The proofs of all results in this section are provided in subsection \ref{subsec_sec_mle}.

\section{The case of several untruncated components}
\label{sec_several_untruncated}

A reviewer raised an interesting question:  Is the exponential family corresponding to a multivariate normal distribution with some, but not all, components truncated also not steep?  To demonstrate that the answer is affirmative, we consider in this section a multivariate normal distribution with one truncated component.  We also note that our proof of non-steepness in this case can be extended, \textit{mutatis mutandis}, to any number of truncated components.

Moreover, it is worthwhile to emphasize that the single-component truncated normal distribution is a multivariate distribution, which is in contrast to the single-component truncated univariate distribution studied in \cite{delCastillo}.  As mentioned earlier, the truncated multivariate normal family is not closed under marginalization \cite{Levine}, so the results established in this section do not follow from the results in \cite{delCastillo,Hegde}.

In this section, we suppose that $X_1$ is truncated to the left at zero and $X_2,\ldots,X_d$ are untruncated.  Similar to Section \ref{sec_non_regularity}, we parametrize the exponential family using the parameter $(\btheta,\bTheta)$.  We partition $\tt$, $\btheta$, and $\bTheta$ as follows: Let $\tt_2 = (t_2,\ldots,t_d)'$ and $\btheta_2 = (\theta_2,\ldots,\theta_d)'$, so that
$$
\tt = \begin{pmatrix} t_1 \\ \tt_2 \end{pmatrix},
\quad
\btheta = \begin{pmatrix} \theta_1 \\ \btheta_2 \end{pmatrix},
\quad
\bTheta = (\Theta_{jk})_{d \times d} =
\begin{pmatrix} \Theta_{11} & \bTheta_{12} \\ \bTheta_{21} & \bTheta_{22} \end{pmatrix},
$$
where $\Theta_{11}$ is $1 \times 1$, $\bTheta_{21} = \bTheta_{12}'$ is $(d-1) \times 1$, and $\bTheta_{22}$ is $(d-1) \times (d-1)$.  Further, we define
\begin{equation}
\label{eq_conditional_thetas}
\theta_{1\cdot 2} := \theta_1 - \bTheta_{12}\bTheta_{22}^{-1}\btheta_2
\ \ {\rm and} \ \
\Theta_{11\cdot 2} := \Theta_{11} - \bTheta_{12}\bTheta_{22}^{-1}\bTheta_{21}.
\end{equation}

Define the Laplace transform of the corresponding exponential family,
\begin{equation}
\label{eq_LT_1_truncated_sec5}
L_T(\btheta,\bTheta) = \int_0^\infty \int_{\bR^{d-1}} \exp(\btheta'\tt - \tt' \bTheta \tt) \, \dd\tt_2 \dd t_1,
\end{equation}
and let $\mcalD = \{(\btheta,\bTheta) \in \bR^d \times \mcalS^{d \times d}: L_T(\btheta,\bTheta) < \infty\}$ be the canonical parameter space associated with \eqref{eq_LT_1_truncated_sec5}.  The following result, which describes $\mcalD$ in detail, is essential for identifying points on $\partial\mcalD$ where the gradient of the corresponding cumulant-generating function remains finite, which in turn implies the non-steepness of the exponential family.

\medskip

\begin{proposition}
\label{prop_space_single}
The single-component truncated multivariate normal distributions constitute a full exponential family with canonical parameter space
$
\mcalD = \mcalD_0 \cup \mcalD_1,
$
where
$$
\mcalD_0 = \{(\btheta,\bTheta) \in \bR^d \times \mcalS^{d \times d}: \theta_{1\cdot 2} < 0, \bTheta_{22} \ \hbox{is positive definite}, \Theta_{11\cdot 2} = 0\}.
$$
and $\mcalD_1 = \bR^d \times \mcalS_d^{d \times d}$.
\end{proposition}

\medskip

We now present the main result of this section, which establishes the non-steepness property of the truncated multivariate normal distributions when one component is truncated.

\medskip

\begin{theorem}
\label{thm_single_comp_non_steep}
The exponential family corresponding to the single-component truncated multivariate normal distributions is not steep.
\end{theorem}

\medskip

The proofs of Proposition \ref{prop_space_single} and Theorem \ref{thm_single_comp_non_steep} are given in subsection \ref{subsec_sec_several_untruncated}.

\section{Concluding remarks}
\label{sec_conclusions}

In this article, we have realized as an exponential family the collection of multivariate normal distributions in which each component is one-sided truncated. We derived the canonical parameter space for this family and showed that the family is full and minimal, but neither regular nor steep.  We also derived the score equations arising from maximum likelihood estimation for the parameters of the truncated multivariate normal distributions and we obtained a necessary condition for the existence of solutions to the score equations.

We note that the articles \cite{Cohen,Manjunath} provided numerical algorithms for solving the score equations; however, those algorithms are based on the tacit assumption that such solutions exist.  In future work, we will farther address the issue of the existence of solutions to the score equations, taking as our guiding light the article by Bedbur, \textit{et al.} \cite{Bedbur_etal}.  We remark that such research necessarily requires new approaches that are able to resolve the situation in which the set of solutions to the score equations contains continuous curves in $\bR^d \times \mcalS^{d \times d}$.

\section{Proofs}
\label{sec_proofs}

\subsection{Proof of Theorem \ref{thm_parameter_space}}
\label{subsec_sec_non_regularity}

Recall that $\mcalS^{d\times d}$ denotes the set of all $d\times d$ symmetric matrices, and $\mcalS_r^{d \times d}$ is the set of all $d \times d$ positive semidefinite (symmetric) matrices of rank $r$, $0 \le r \le d$.  Let us define $\mcalB = \{(\btheta,\bTheta) \in \bR^d \times \mcalS^{d \times d}\} \equiv \bR^d \times \mcalS^{d \times d}$ and $\mcalB_r = \{(\btheta,\bTheta) \in \bR^d \times \mcalS_r^{d \times d}\} \equiv \bR^d \times \mcalS_r^{d \times d}$.  

Suppose that $\bTheta \in \mcalS_0^{d \times d}$, i.e., $\bTheta = \bzero$.  Then the integral \eqref{eq_Laplace_transform} converges if and only if $\btheta < \bzero$.  Therefore on defining
\begin{equation}
\label{eq_A_0}
\mcalD_0 := \{(\btheta,\bTheta) \in \bR^d \times \mcalS_0^{d \times d}: \btheta < \bzero\} \equiv (-\infty,0)^d \times \{\bzero\},
\end{equation}
we find that $\mcalD_0 = \mcalB_0 \cap \mcalD \subset \mcalD$.  Moreover $\mcalD_0$ is not open in $\bR^d \times \mcalS^{d \times d}$ since its complement is not closed.

Suppose that $\bTheta \in \mcalS_r^{d \times d}$, where $1 \le r \le d-1$.  Then there exists a $d \times d$ orthogonal matrix $\HH$ such that $\HH'\bTheta\HH$ is diagonal and of rank $r$.  We write this in the form $\HH'\bTheta\HH = \diag(\xi_1,\ldots,\xi_r,0,\ldots,0)$, where $\xi_1,\ldots,\xi_r > 0$.  Making the change-of-variables $\ww \to \HH\ww$ in \eqref{eq_Laplace_transform}, we obtain
\begin{equation}
\label{eq_LT_orthogonal}
L_T(\btheta,\bTheta) = \int_{\HH\ww > \bzero} \exp(\btheta'\HH\ww - \ww'\HH'\bTheta\HH\ww) \, \dd\ww.
\end{equation}
Let $\bXi^{-1} = 2 \diag(\xi_1,\ldots,\xi_r)$, and partition $\ww$ into sub-vectors $\ww_1$ and $\ww_2$ such that $\ww' = (\ww_1',\ww_2')$, where $\ww_1$ is $r \times 1$ and $\ww_2$ is $(d-r) \times 1$.  Also let $\HH_1$ and $\HH_2$ be the first $r$ and the last $d-r$ columns, respectively, of $\HH$.  Then
$$
\HH = (\HH_1 \vdots \hspace{1pt} \HH_2), \quad \ww'\HH'\bTheta\HH\ww = \tfrac12 \ww_1'\bXi^{-1}\ww_1, \quad \HH\ww = \HH_1\ww_1 + \HH_2\ww_2,
$$
and it follows from \eqref{eq_LT_orthogonal} that
$$
L_T(\btheta,\bTheta) = \int_{\HH_1\ww_1 + \HH_2\ww_2 > \bzero} \exp(\btheta'\HH_1\ww_1 - \tfrac12 \ww_1'\bXi^{-1}\ww_1 + \btheta'\HH_2\ww_2) \dd\ww.
$$
By the standard device of ``completing the square'' of the quadratic form, we obtain
$$
\btheta'\HH_1\ww_1 - \tfrac12 \ww_1'\bXi^{-1}\ww_{1} \equiv -\tfrac12(\ww_1 - \bXi\HH_1'\btheta)'\bXi^{-1}(\ww_1 - \bXi\HH_1'\btheta) + \tfrac12 \btheta'\HH_1\bXi\HH_1'\btheta.
$$
Let $n(\ww_1;\bXi\HH_1'\btheta,\bXi)$, $\ww_1 \in \bR^r$, denote the probability density function of the multivariate normal random vector $\WW_1 \sim N_r(\bXi\HH_1'\btheta,\bXi)$; then 
\begin{align*}
&\exp(-\tfrac12 \btheta'\HH_1\bXi\HH_1'\btheta) L_T(\btheta,\bTheta) \\
&= \int_{\HH_1\ww_1 + \HH_2\ww_2 > \bzero} \exp[-\tfrac12(\ww_1 - \bXi\HH_1'\btheta)'\bXi^{-1}(\ww_1 - \bXi\HH_1'\btheta)] \exp(\btheta'\HH_2\ww_2) \dd\ww \\
&\propto \int_{\HH_1\ww_1 + \HH_2\ww_2 > \bzero} n(\ww_1;\bXi\HH_1'\btheta,\bXi) \exp(\btheta'\HH_2\ww_2) \dd\ww_1 \dd\ww_2.
\end{align*}
Also denoting by $\chi(A)$ the indicator function of any interval $A \subset \bR$, we obtain
\begin{align*}
&\exp(-\tfrac12 \btheta'\HH_1\bXi\HH_1'\btheta) L_T(\btheta,\bTheta) \\
\ &\propto \int_{\bR^r} n(\ww_1;\bXi\HH_1'\btheta,\bXi) \bigg[\int_{\bR^{d-r}} \chi(\HH_2\ww_2 > -\HH_1\ww_1) \exp(\btheta'\HH_2\ww_2) \dd\ww_2\bigg] \dd\ww_1.
\end{align*}
It follows that if $\bTheta \in \mcalS_r^{d \times d}$ then $(\btheta,\bTheta) \in \mcalD$ if and only if $\btheta$ satisfies the condition
$$
\E_{\WW_1} \int_{\bR^{d-r}} \chi(\HH_2\ww_2 > -\HH_1\WW_1) \exp(\btheta'\HH_2\ww_2) \dd\ww_2 < \infty.
$$
Consequently, defining
\begin{multline}
\label{eq_A_r}
\mcalD_r := \bigg\{(\btheta,\bTheta) {\hskip-0.5pt} \in \bR^d {\hskip-0.5pt} \times {\hskip-0.5pt} \mcalS_r^{d \times d}: \\
\E_{\WW_1} {\hskip-0.44pt} \int_{\bR^{d-r}} \chi(\HH_2\ww_2 > -\HH_1\WW_1) \exp(\btheta'\HH_2\ww_2) \dd\ww_2 < \infty\bigg\},
\end{multline}
where $1 \le r \le d-1$, then $\mcalD_r = \mcalB_r \cap \mcalD \subset \mcalD$.

To show that $\mcalS_r^{d \times d}$ is not open in $\mcalS^{d \times d}$, suppose that $\bTheta \in \mcalS_r^{d \times d}$ with $1 \le r \le d-1$.  For every $\epsilon > 0$ and $\ww \in \bR^d$ such that $\ww \neq \bzero$, there holds 
$$
\ww'(\bTheta + \epsilon \II_d)\ww > 0.
$$
Hence, in every neighborhood of $\bTheta$ there exists a positive definite matrix, $\bTheta + \epsilon \II_d$, that is not in $\mcalS_r^{d \times d}$.  Therefore $\mcalS_r^{d \times d}$ is not open in $\mcalS^{d \times d}$, and it follows that $\mcalD_r$, for $1 \le r \le d-1$, is not open.  

Next, consider the case $r=d$.  Then $\bTheta \in \mcalS_d^{d \times d}$, and the integral \eqref{eq_Laplace_transform} is a constant multiple of a multivariate normal probability so the integral converges for all $\btheta \in \bR^d$.  Therefore, defining 
\begin{equation}
\label{eq_A_d}
\mcalD_d := \bR^d \times \mcalS_d^{d \times d},
\end{equation}
we obtain $\mcalD_d = \mcalB_d \cap \mcalD$, an open subset of $\mcalD$.  

Finally, consider the case in which $\bTheta \not\in \bigcup_{r=0}^d \mcalS_r^{d \times d}$, i.e., $\bTheta$ is not positive semidefinite.  Then $\bTheta$ has at least one negative eigenvalue.  Suppose that the eigenvalues of $\bTheta$ are $-\xi_1,\ldots,-\xi_r, \xi_{r+1},\ldots,\xi_d$ where $r \ge 1$, $\xi_1,\ldots,\xi_r > 0$, and $\xi_{r+1},\ldots,\xi_d \ge 0$.  Let $\bXi_1 = \diag(\xi_1,\ldots,\xi_r)$, $\bXi_2 = \diag(\xi_{r+1},\ldots,\xi_d)$, and apply an orthogonal transformation, $\ww \to \HH\ww$, to rewrite the Laplace transform \eqref{eq_Laplace_transform} as in \eqref{eq_LT_orthogonal}.  Also partition $\ww$ and $\HH$ as before, with
$$
\ww = \begin{pmatrix}\ww_1 \\ \ww_2\end{pmatrix}, \qquad \HH = (\HH_1 \vdots \hspace{1pt} \HH_2),
$$
where $\ww_1$ is $r \times 1$, $\ww_2$ is $(d-r) \times 1$, $\HH_1$ is $d \times r$, and $\HH_2$ is $d \times (d-r)$.  Then
\begin{align*}
&L_T(\btheta,\bTheta) \\
&= \int_{\HH_1\ww_1 + \HH_2\ww_2 > \bzero} \exp(\btheta'\HH_1\ww_1 + \ww_1'\bXi_1\ww_1 + \btheta'\HH_2\ww_2 - \ww_2'\bXi_2\ww_2) \dd\ww_1 \dd\ww_2.
\end{align*}
Since
$$
\{(\ww_1,\ww_2): \HH_1\ww_1 + \HH_2\ww_2 > \bzero\} \supseteq \{(\ww_1,\ww_2): \HH_1\ww_1 > \bzero, \HH_2\ww_2 > \bzero\}
$$
then
\begin{align*}
&L_T(\btheta,\bTheta) \\
&\ge \int_{\HH_1\ww_1 > \bzero, \HH_2\ww_2 > \bzero} \exp(\btheta'\HH_1\ww_1 + \ww_1'\bXi_1\ww_1 + \btheta'\HH_2\ww_2 - \ww_2'\bXi_2\ww_2) \dd\ww_1 \dd\ww_2 \\
&= \int_{\HH_1\ww_1 > \bzero} \exp(\btheta'\HH_1\ww_1 + \ww_1'\bXi_1\ww_1) \dd\ww_1 \\
& \qquad \cdot \int_{\HH_2\ww_2 > \bzero} \exp(\btheta'\HH_2\ww_2 - \ww_2'\bXi_2\ww_2) \dd\ww_2.
\end{align*}
Since $\bXi_1$ is positive definite then $\btheta'\HH_1\ww_1 + \ww_1'\bXi_1\ww_1 \to \infty$ as $\|\ww_1\| \to \infty$.  Therefore the first integral with respect to $\ww_1$ in the last expression diverges, so the Laplace transform $L_T(\btheta,\bTheta)$ also diverges.  Hence,
$$
\bigg\{(\btheta,\bTheta) \in \bR^d \times \mcalS^{d \times d}: \bTheta \not\in \bigcup_{r=0}^d \mcalS_r^{d \times d}\bigg\} \subset \mcalD^{c},
$$
equivalently, $(\mcalB\setminus\bigcup_{r=0}^d\mcalB_r) \bigcap\mcalD=\emptyset$.  

Note that the collection $\{\mcalB_0,\ldots,\mcalB_d,\mcalB \setminus (\bigcup_{r=0}^d\mcalB_r)\}$ forms a partition of $\mcalB$.  Moreover, observe that $\mcalD \subseteq \mcalB$.  By the standard partition-decomposition identity, we obtain 
$$
\mcalD = \bigcup_{r=0}^d(\mcalB_r \cap \mcalD) \, \bigcup \, \Big(\big(\mcalB \setminus \bigcup_{r=0}^d \, \mcalB_r\big) \cap \mcalD\Big) = \bigcup_{r=0}^d(\mcalB_r\cap \mcalD) = \bigcup_{r=0}^d\mcalD_r,
$$
as asserted in \eqref{eq_D_nat_param_space}.  

Since each $\mcalD_r$, $0 \le r \le d-1$, is not open and $\mcalD_0,\ldots,\mcalD_d$ are mutually disjoint then their union $\mcalD$ is not open.  Therefore the family $\mcalP_\mcalD$ is not regular.

On the other hand, $\mcalP_\mcalD$ is full since $\mcalD$ is of dimension $k = d+\tfrac12 d(d+1)$, which is the dimension of the statistic $T$ in \eqref{eq_T_of_tt}.  Moreover, $\mcalP_\mcalD$ is minimal since $k$ is also the dimension of the canonical statistic $T(\tt)$.  

Finally, the formula \eqref{eq_trunc_normal_nef_pdf} for the density function $p(\tt;\btheta,\bTheta)$, when $(\btheta,\bTheta) \in \mcalD_d$, follows from \eqref{eq_KT_cgf}.
$\qed$

\subsection{Proof of Theorem \ref{thm_non_steep}}
\label{subsec_sec_non_steep}

Throughout the proof of the next result, we will sometimes interchanging integrals and derivatives.  In all such instances, the interchange can be validated by applying standard exponential family results from Barndorff-Nielsen \cite{Barndorff} or classical mathematical analysis from Burkill and Burkill \cite[p.~290]{Burkill}.  Further, all interchanges of limits and integrals can be justified by a straightforward appeal to the dominated convergence theorem.  

\medskip

\noindent\textit{Proof of Theorem \ref{thm_non_steep}}.  
Let $\btheta_0 \in \bR^d$ be a fixed vector such that $\btheta_0 < \bzero$, and let $\bomega_0 = (\btheta_0,\bzero)$.  It is elementary from \eqref{eq_Laplace_transform} and \eqref{eq_D_space} that $L_T(\bomega_0) < \infty$, therefore $\bomega_0 \in \mcalD$.  Now consider the sequence $\bomega_n = (\btheta_0,n^{-1} I_d)$, $n=1,2,3,\ldots$.  Since the matrix $n^{-1} \II_d$ is positive definite then it again follows from \eqref{eq_Laplace_transform} and \eqref{eq_D_space} that $L_T(\bomega_n) < \infty$.  Hence, by \eqref{eq_A_d}, $\bomega_n \in \mcalD_d \equiv \hbox{Int}(\mcalD_d)$, since $\mcalD_d$ is open.  Since $\bomega_n \to \bomega_0$ as $n \to \infty$ then it follows that $\bomega_0 \in \cl(\mcalD_d) \subseteq \cl(\Omega)$.

To show that $\bomega_0 \in \cl(\Omega^c)$, consider a new sequence $\bomega_n = (\btheta_0,-n^{-1} \II_d)$, $n=1,2,3,\ldots$.  Since $-n^{-1} \II_d$ is negative definite then it follows from \eqref{eq_Laplace_transform} that $L_T(\bomega_n) = \infty$.  By \eqref{eq_D_space}, $\bomega_n \in \mcalD^c \subseteq \Omega^c$, and since $\bomega_n \to \bomega_0$ as $n \to \infty$ then $\bomega_0 \in \cl(\Omega^c)$.  Therefore $\bomega_0 \in \cl(\Omega) \cap \cl(\Omega^c) \equiv \partial\Omega$.

Next, we apply \eqref{eq_steep} to show that the truncated multivariate normal family is not steep at $\bomega_0$.  Starting with \eqref{eq_KT_cgf}, we have 
$$
\nabla_\btheta K_T(\btheta,\bTheta) = \frac{\int_{\tt > \bzero} \tt \exp(\btheta'\tt - \tt' \bTheta \tt) \, \dd\tt}{\int_{\tt > \bzero} \exp(\btheta'\tt - \tt' \bTheta \tt) \, \dd\tt}
$$
and
$$
\nabla_\bTheta K_T(\btheta,\bTheta) = - \frac{\int_{\tt > \bzero} \tt\tt' \exp(\btheta'\tt - \tt' \bTheta \tt) \, \dd\tt}{\int_{\tt > \bzero} \exp(\btheta'\tt - \tt' \bTheta \tt) \, \dd\tt}.
$$
Consider a sequence $(\btheta_n,\bTheta_n) \in \bR^d \times \mcalS_d^{d \times d}$, $n=1,2,3,\ldots$ where $\btheta_n = \btheta_0$ for all $n$ and $\bTheta_n \to \bzero$ as $n \to \infty$.  It is straightforward that $(\btheta_n,\bTheta_n) \in \mcalD$ for all $n$.

Further,
$$
\nabla_\btheta K_T(\btheta_n,\bTheta_n) := \nabla_\btheta K_T(\btheta,\bTheta){\Big|}_{\btheta=\btheta_0,\bTheta = \bTheta_n} = \dfrac{\int_{\tt > \bzero} \tt \exp(\btheta_0'\tt - \tt'\bTheta_n\tt) \, \dd\tt}{\int_{\tt > \bzero} \exp(\btheta_0'\tt - \tt'\bTheta_n\tt) \, \dd\tt},
$$
and the numerator and denominator integrals both exist since $\bTheta_n$ is positive definite.  Since $\bTheta_n \to \bzero$ as $n \to \infty$, then,
$$
\nabla_\btheta K_T(\btheta_n,\bTheta_n) \to \dfrac{\int_{\tt > \bzero} \tt \exp(\btheta_0'\tt) \, \dd\tt}{\int_{\tt > \bzero} \exp(\btheta_0'\tt) \, \dd\tt},
$$
and these integrals also are finite since $\btheta_0 < \bzero$.  Denote by $\theta_{0,j}$ the $j$th component of $\btheta_0$; on calculating the latter integrals, it is seen that $\nabla_\btheta K_T(\btheta_n,\bTheta_n)$ converges to a vector whose $j$th component is
$$
\dfrac{\int_0^\infty t_j \exp(\theta_{0,j} t_j) \, \dd t_j}{\int_0^\infty  \exp(\theta_{0,j} t_j) \, \dd t_j} = \frac{(-\theta_{0,j})^{-2}}{(-\theta_{0,j})^{-1}} = - \theta_{0,j}^{-1},
$$
$j=1,\ldots,d$, and therefore
\begin{equation}
\label{eq_nabla_theta_KT}
\|\nabla_\btheta K_T(\btheta_n,\bTheta_n)\|^2 \to \sum_{j=1}^d \theta_{0,j}^{-2}
\end{equation}
as $n \to \infty$.  Note that we also have
$$
\btheta_0' \nabla_\btheta K_T(\btheta_n,\bTheta_n) \to \sum_{j=1}^d \theta_{0,j} (- \theta_{0,j}^{-1}) = -d
$$
as $n \to \infty$, which generalizes a result of del Castillo \cite[p.~62]{delCastillo}.

Next,
\begin{align*}
\nabla_\bTheta K_T(\btheta_n,\bTheta_n) &:= \nabla_\bTheta K_T(\btheta,\bTheta)\Big|_{\btheta=\btheta_0,\bTheta = \bTheta_n} \\
&= - \dfrac{\int_{\tt > \bzero} \tt\tt' \exp(\btheta_0'\tt - \tt'\bTheta_n\tt) \, \dd\tt}{\int_{\tt > \bzero} \exp(\btheta_0'\tt - \tt'\bTheta_n\tt) \, \dd\tt}.
\end{align*}
Therefore, as $n \to \infty$,
$$
\nabla_\bTheta K_T(\btheta_n,\bTheta_n) \to - \dfrac{\int_{\tt > \bzero} \tt\tt' \exp(\btheta_0'\tt) \, \dd\tt}{\int_{\tt > \bzero} \exp(\btheta_0'\tt) \, \dd\tt},
$$
and these integrals are finite since $\btheta_0 < \bzero$.  On calculating these integrals we find that $\nabla_\bTheta K_T(\btheta_n,\bTheta_n)$ converges to $\MM$, a $d \times d$ matrix whose $(i,j)$th entry is
$$
m_{i,j} = - \dfrac{\int_0^\infty \int_0^\infty t_i t_j \exp(\theta_{0,i} t_i + \theta_{0,j} t_j) \, \dd t_i \dd t_j}{\int_0^\infty \int_0^\infty \exp(\theta_{0,i} t_i + \theta_{0,j} t_j) \, \dd t_i \dd t_j}
= \theta_{0,i}^{-1} \theta_{0,j}^{-1},
$$
$i,j=1,\ldots,d$.  Hence, as $n \to \infty$,
\begin{equation}
\label{eq_nabla_Theta_KT}
\|\nabla_\bTheta K_T(\btheta_n,\bTheta_n)\|_F^2 \to \|\MM\|_F^2 = \sum_{i,j=1}^d m_{i,j}^2 = \sum_{i,j=1}^d \theta_{0,i}^{-2} \theta_{0,j}^{-2} = \Big(\sum_{j=1}^d \theta_{0,j}^{-2}\Big)^2.
\end{equation}
Consequently it follows from \eqref{eq_nabla_theta_KT} and \eqref{eq_nabla_Theta_KT} that, for $\btheta < \bzero$,
$$
\|\nabla_\bomega K_T(\btheta_n,\bTheta_n)\|^2 \to \sum_{j=1}^d \theta_{0,j}^{-2} + \Big(\sum_{j=1}^d \theta_{0,j}^{-2}\Big)^2,
$$
a finite limit, as $n \to \infty$.
$\qed$

\subsection{Proofs for Section \ref{sec_mle}}
\label{subsec_sec_mle}

\noindent\textit{Proof of Theorem \ref{thm_existence_mles}}. 
(i) In the case of the untruncated multivariate normal distributions, it is well known that the pair $(\TT_1,\TT_2)$ is a sufficient statistic.  Conditioning the untruncated normal distributions, or any exponential family, on a subset of the sample space does not affect this sufficiency property since only the normalizing constant in the density needs to be modified.  Therefore $(\TT_1,\TT_2)$ is also sufficient in the truncated case.  

Alternatively, one may proceed directly by showing, as in the classical untruncated case, that the log-likelihood function is 
\begin{align*}
\ell(\bmu,\bSigma) &= \sum_{j=1}^n \log f(\xx_j;\bmu,\bSigma) \nonumber \\
&= -n \log C(\bmu,\bSigma) - \frac{1}{2}n \tr \bSigma^{-1}(\TT_2 - \TT_1 \bmu' - \bmu \TT_1' + \bmu\bmu'),
\end{align*}
which shows that $(\TT_1,\TT_2)$ is sufficient for $(\bmu,\bSigma)$.  This proves (i).

(ii) For any random sample $\xx_1,\ldots,\xx_n$, it is a classical identity that  
\begin{equation}
\label{eq_T2T1}
\TT_2 - \TT_1 \TT_1' = n^{-1} \sum_{j=1}^n \xx_j \xx_j' - \bar{\xx}\bar{\xx}' \equiv n^{-1} \sum_{j=1}^n (\xx_j - \bar{\xx})(\xx_j - \bar{\xx})',
\end{equation}
which evidently is positive semidefinite.  Therefore $(\TT_1,\TT_2) \in \mcalC_d$.  

To show that $\mcalC_d$ is closed, let $(\btheta_n,\bTheta_n)$, $n=1,2,3,\ldots$, be a sequence in $\mcalC_d$ that converges to $(\btheta,\bTheta)$, i.e., $\btheta_n \to \btheta$ and $\bTheta_n \to \bTheta$.  Since $(\btheta_n,\bTheta_n) \in \mcalC_d$ then $\bTheta_n - \btheta_n\btheta_n' \succcurlyeq \bzero$ for all $n$.  Since the set of all positive semidefinite $d \times d$ matrices is closed then $\bTheta_n - \btheta_n\btheta_n' \to \bTheta - \btheta\btheta' \succcurlyeq \bzero$.  Therefore $(\btheta,\bTheta) \in \mcalC_d$, hence $\mcalC_d$ is closed.  

Next, for $\tau \in (0,1)$ and $(\bmu_1,\bSigma_1), (\bmu_2,\bSigma_2) \in \mcalC_d$, define the convex combination 
\begin{align*}
(\bmu_3,\bSigma_3) &= \tau (\bmu_1,\bSigma_1) + (1-\tau)(\bmu_2,\bSigma_2) \\
&= (\tau \bmu_1 + (1-\tau) \bmu_2,\tau \bSigma_1 + (1-\tau)\bSigma_2).
\end{align*}
A straightforward algebraic calculation reveals that
\begin{align*}
\bSigma_3 - \bmu_3\bmu_3' &\equiv \tau \bSigma_1 + (1-\tau)\bSigma_2 - (\tau \bmu_1 + (1-\tau) \bmu_2)(\tau \bmu_1 + (1-\tau) \bmu_2)' \\
&= \tau(\bSigma_1 - \bmu_1\bmu_1') + (1-\tau)(\bSigma_2 - \bmu_2\bmu_2') \\
& \qquad\qquad\qquad\qquad\qquad + \tau(1-\tau)(\bmu_1 - \bmu_2)(\bmu_1 - \bmu_2)'.
\end{align*}
This shows that $\bSigma_3 - \bmu_3\bmu_3'$ is positive semidefinite since it is a sum of three positive semidefinite matrices.  Therefore $(\bmu_3,\bSigma_3) \in \mcalC_d$, so $\mcalC_d$ is convex.

(iii) Observe that the random sample $(\xx_1,\ldots,\xx_n)$ has a joint probability density on the underlying Euclidean space $\bR^{dn}$.  Therefore by a result of Malley \cite[p.~344]{Malley}, the probability distribution of $(\xx_1,\ldots,\xx_n)$ assigns zero probability to the zeros of any non-trivial polynomial in the components of $(\xx_1,\ldots,\xx_n)$.  For $n-1 \ge d$, equivalently, $n > d$, it follows from \eqref{eq_T2T1} that $\TT_2 - \TT_1\TT_1'$ is of full rank, almost surely.  Applying Malley's result to $\det(\TT_2 - \TT_1\TT_1')$, which is a non-trivial polynomial in the components of $(\xx_1,\ldots,\xx_n)$, it follows that $\P\big(\det(\TT_2 - \TT_1\TT_1') = 0\big) = 0$, equivalently $\det(\TT_2 - \TT_1\TT_1') > 0$, so $\TT_2 - \TT_1\TT_1'$ is positive definite, almost surely.  

Recall from Theorem \ref{thm_parameter_space} that the family $\mcalP_\mcalD$ is minimal and full.  Therefore the log-likelihood function, as a function of the canonical parameter $(\btheta,\bTheta)$, is strictly concave on $\mcalD$, the convex natural parameter space.  Hence, if there exists on $\mcalD$ a maximizer of $\ell(\btheta,\bTheta)$ then the maximizer is unique \cite[p.~151, Theorem 9.13]{Barndorff}.  By making the bijective transformation from $(\btheta,\bTheta)$ to $(\bmu,\bSigma)$, given by $\bmu = (2\Theta)^{-1}\btheta$ and $\bSigma = (2\Theta)^{-1}$, we find that if there exists a maximizer of the function $\ell(\bmu,\bSigma)$ then that maximizer also is unique.  
$\qed$

\bigskip

\noindent\textit{Proof of Lemma \ref{lem_mean_covariance}}.
For $\tt = (t_1,\ldots,t_d)' \in \bR^d$, it follows from \eqref{eq_trunc_normal_pdf} that the moment-generating function of $\XX$ is
$$
\E(e^{\tt'\XX}) = \frac{1}{[C(\bmu,\bSigma)]} \int_{\bR^d_{+}} \exp[\tt'\xx - \tfrac12(\xx-\bmu)'\bSigma^{-1}(\xx-\bmu)] \dd \xx.
$$
Applying the usual approach of completing the square, we obtain the algebraic identity,
$$
\tt'\xx - \tfrac12(\xx-\bmu)'\bSigma^{-1}(\xx-\bmu) \equiv \tt'\bmu + \tfrac12 \tt' \bSigma \tt -\tfrac12(\xx-\bmu-\bSigma \tt)'\bSigma^{-1}(\xx-\bmu-\bSigma \tt).
$$
Hence
\begin{align*}
\E(e^{\tt'\XX}) &= \frac{1}{[C(\bmu,\bSigma)]} \exp(\tt'\bmu + \tfrac12 \tt' \bSigma \tt) \\
& \qquad \times \int_{\bR^d_{+}} \exp[-\tfrac12(\xx-\bmu-\bSigma \tt)'\bSigma^{-1}(\xx-\bmu-\bSigma \tt)] \dd \xx \\
&= \exp(\tt'\bmu + \tfrac12 \tt' \bSigma \tt) \, \frac{C(\bmu + \bSigma \tt,\bSigma)}{[C(\bmu,\bSigma)]},
\end{align*}
and therefore
\begin{equation}
\label{eq_cgf}
\log \e(e^{\tt'\XX}) = K(\tt) = \tt'\bmu + \frac{1}{2} \tt' \bSigma \tt + \log C(\bmu + \bSigma \tt,\bSigma) - \log C(\bmu,\bSigma).
\end{equation}
By the chain rule,
$$
\nabla_{\tt} [C(\bmu + \bSigma \tt,\bSigma)] = \bSigma \, (\nabla C)(\bmu + \bSigma \tt,\bSigma);
$$
therefore, by differentiating \eqref{eq_cgf}, we obtain
\begin{align}
\label{eq_nabla_K}
\frac{\E(e^{\tt'\XX} \XX)}{\E(e^{\tt'\XX})} = \nabla_{\tt} K(\tt) &= \bmu + \bSigma \tt + \frac{1}{C(\bmu + \bSigma \tt,\bSigma)} \ \bSigma \, (\nabla C)(\bmu + \bSigma \tt,\bSigma) \nonumber \\
&= \bmu + \bSigma \, \Big(\tt + \frac{1}{C(\bmu + \bSigma \tt,\bSigma)} \, (\nabla C)(\bmu + \bSigma \tt,\bSigma)\Big).
\end{align}
Setting $\tt = \bzero$ we obtain
\begin{equation}
\label{eq_nu}
\bnu = \E(\XX) = \nabla K(\tt)\big|_{\tt=\boldsymbol0} = \bmu + \frac{1}{C(\bmu,\bSigma)} \ \bSigma \, (\nabla C)(\bmu,\bSigma),
\end{equation}
and this establishes \eqref{eq_nu_ML}.

Next, it follows from \eqref{eq_nabla_K} that
\begin{align}
\label{eq_2nd_deriv_K}
\nabla_{\tt} \nabla_{\tt}' K(\tt) &\equiv \nabla_{\tt} [\nabla_{\tt} K(\tt)]' \nonumber \\
&= \nabla_{\tt} \Big[\bmu + \bSigma \, \Big(\tt + \frac{1}{C(\bmu + \bSigma \tt,\bSigma)} \, (\nabla C)(\bmu + \bSigma \tt,\bSigma)\Big)\Big]' \nonumber \\
&= \nabla_{\tt} \Big(\tt' + \frac{1}{C(\bmu + \bSigma \tt,\bSigma)} \, (\nabla C)'(\bmu + \bSigma \tt,\bSigma)\Big) \bSigma.
\end{align}
Since $\nabla_{\tt} \tt' = \II_d$ then, by the chain rule,
\begin{align}
\label{eq_2nd_deriv_K_2}
&\nabla_{\tt}\Big(\frac{1}{C(\bmu + \bSigma \tt,\bSigma)} \, (\nabla C)'(\bmu + \bSigma \tt,\bSigma)\Big) \nonumber \\
&= \frac{1}{[C(\bmu + \bSigma \tt,\bSigma)]^2} \, \bSigma \big[C(\bmu + \bSigma \tt,\bSigma) \, (\nabla \nabla' C)(\bmu + \bSigma \tt,\bSigma) \nonumber \\
& \qquad\qquad\qquad\qquad\qquad\qquad - (\nabla C)(\bmu + \bSigma \tt,\bSigma) \, (\nabla C)'(\bmu + \bSigma \tt,\bSigma)\big].
\end{align}
Substituting \eqref{eq_2nd_deriv_K_2} into \eqref{eq_2nd_deriv_K} and evaluating the resulting expression at $\tt = \bzero$, we obtain
\begin{align}
\label{eq_Lambda_2}
\bLambda &= \nabla_{\tt} \nabla_{\tt}' K(\tt)\Big|_{\tt = \bzero} \nonumber \\
&= \bSigma + [C(\bmu,\bSigma)]^{-2} \bSigma \big[C(\bmu,\bSigma)  (\nabla \nabla' C)(\bmu,\bSigma) \nonumber \\
& \qquad\qquad\qquad\qquad\qquad\qquad - (\nabla C)(\bmu,\bSigma)  (\nabla C)'(\bmu,\bSigma) \big]\bSigma \nonumber \\
&= \bSigma + [C(\bmu,\bSigma)]^{-1} \bSigma (\nabla \nabla' C)(\bmu,\bSigma) \bSigma \nonumber \\
& \qquad\qquad\qquad- [C(\bmu,\bSigma)]^{-2} \bSigma (\nabla C)(\bmu,\bSigma) \cdot (\nabla C)'(\bmu,\bSigma) \bSigma.
\end{align}
By \eqref{eq_nu_ML},
\begin{align*}
[C(\bmu,\bSigma)]^{-2} & \bSigma (\nabla C)(\bmu,\bSigma) \cdot (\nabla C)'(\bmu,\bSigma) \bSigma \\
&\equiv [C(\bmu,\bSigma)]^{-1} \bSigma (\nabla C)(\bmu,\bSigma) \cdot \big([C(\bmu,\bSigma)]^{-1} \bSigma (\nabla C)(\bmu,\bSigma)]\big)' \\
&= (\bnu - \bmu)(\bnu - \bmu)',
\end{align*}
and by substituting this expression into \eqref{eq_Lambda_2}, we obtain \eqref{eq_Lambda_ML}.
$\qed$

\bigskip

\noindent\textit{Proof of Theorem \ref{thm_score_equations}}.
In deriving the score equations, we treat the vectors $\xx_1,\ldots,\xx_n$ as fixed at their sample values, and $\bmu$ and $\bSigma$ are viewed temporarily as variables.  Then the log-likelihood function corresponding to the random sample $\xx_1,\ldots,\xx_n$ is
\begin{align}
\label{eq_llh_function2}
\ell(\bmu,\bSigma) &= \sum_{j=1}^n \log f(\xx_j;\bmu,\bSigma) \nonumber \\
&= -n \log C(\bmu,\bSigma) - \frac{1}{2} \sum_{j=1}^n (\xx_j-\bmu)'\bSigma^{-1}(\xx_j-\bmu),
\end{align}
and therefore
\begin{equation}
\label{eq_mu_score}
\nabla_\bmu \ell(\bmu,\bSigma) = - n \frac{\nabla_\bmu C(\bmu,\bSigma)}{C(\bmu,\bSigma)} + n \, \bSigma^{-1}(\overline{\xx}-\bmu).
\end{equation}
By \eqref{eq_nu},
$$
\frac{\nabla_\bmu C(\bmu,\bSigma)}{C(\bmu,\bSigma)} = \bSigma^{-1} (\bnu - \bmu),
$$
and by substituting this result into \eqref{eq_mu_score} we obtain
\begin{align*}
\nabla_\bmu \ell(\bmu,\bSigma) &= -n\, \bSigma^{-1}(\bnu - \bmu) + n \,\bSigma^{-1}(\overline{\xx} - \bmu) \\
&= -n\, \bSigma^{-1} (\bnu - \overline{\xx}).
\end{align*}

As $(\widehat\bmu,\widehat\bSigma)$ is a stationary point of $\ell(\bmu,\bSigma)$, we obtain
$$
-n\, \widehat\bSigma^{-1} (\widehat\bnu - \overline{\xx}) = \nabla_\bmu \ell(\bmu,\bSigma)\Big|_{(\bmu,\bSigma)=(\widehat\bmu,\widehat\bSigma)} = \bzero,
$$
hence $\widehat\bnu = \overline{\xx}$.  Substituting for $\widehat\bnu$ from \eqref{eq_nu_MLE}, we obtain \eqref{eq_score_eqs_nu}.

Next, we derive the score equations for $\bSigma$ in terms of $\bPsi = \bSigma^{-1}$.  Denoting $\ell(\bmu,\bSigma)$ by $\tilde\ell(\bmu,\bPsi)$, $C(\bmu,\bSigma)$ by $\tilde{C}(\bmu,\bPsi)$, and using the definition of $\SS(\bmu)$ from \eqref{eq_def_A}, we find that the log-likelihood function \eqref{eq_llh_function2} equals
\begin{equation}
\label{eq_score_funct_Psi}
\tilde\ell(\bmu,\bPsi) = -n \log \tilde{C}(\bmu,\bPsi) - \frac{1}{2} n \tr \bPsi \SS(\bmu).
\end{equation}
Denote by $\psi_{i,j}$ the $(i,j)$th entry of $\bPsi$, $i,j=1,\ldots,d$.  Since
$$
\tilde{C}(\bmu,\bPsi) = \int_{\ww \in \mathbb{R}^d_+} \exp\big[-\tfrac12 (\ww-\bmu)'\bPsi(\ww-\bmu)\big] \dd \ww
$$
then, by differentiating \eqref{eq_score_funct_Psi}, we obtain
\begin{equation}
\label{eq_score_Psi}
\frac{\partial}{\partial\psi_{i,j}} \tilde\ell(\bmu,\bPsi) = - \frac{n}{\tilde{C}(\bmu,\bPsi)} \frac{\partial}{\partial\psi_{i,j}} \tilde{C}(\bmu,\bPsi) - \frac{n}{2} \frac{\partial}{\partial\psi_{i,j}} \tr \bPsi \SS(\bmu).
\end{equation}
For $i,j=1,\ldots,d$, denote by $a_{i,j}$ the $(i,j)$th element of $\SS(\bmu)$.  Since $\SS(\bmu)$ is symmetric then
\begin{equation}
\label{eq_deriv_trace}
\frac{\partial}{\partial\psi_{i,j}} \tr \bPsi \SS(\bmu) = (2-\delta_{i,j}) \,a_{i,j}.
\end{equation}
Differentiating under the integral sign we obtain
\begin{align}
\label{eq_deriv_C_mu_Psi}
&\frac{\partial}{\partial\psi_{i,j}} \tilde{C}(\bmu,\bPsi) \nonumber \\
&= \int_{\ww \in \mathbb{R}^d_+} \frac{\partial}{\partial\psi_{i,j}} \exp\big[-\tfrac12(\ww-\bmu)'\bPsi(\ww-\bmu)\big] \dd\ww \nonumber \\
&= -\frac12 \int_{\ww \in \mathbb{R}^d_+} (2-\delta_{i,j})(w_i-\mu_i)(w_j-\mu_j) \exp{\hskip-0.35pt}\big[{\hskip-2pt}-{\hskip-1pt}\tfrac12 (\ww-\bmu)'\bPsi(\ww-\bmu)\big] \dd\ww \nonumber \\
&= -\frac12\, (2-\delta_{i,j})\, \tilde{C}(\bmu,\bPsi)\, \E[(X_i-\mu_i)(X_j-\mu_j)].
\end{align}
Substituting \eqref{eq_deriv_trace} and \eqref{eq_deriv_C_mu_Psi} into \eqref{eq_score_Psi}, we obtain
$$
\frac{\partial}{\partial\psi_{i,j}}\tilde\ell(\bmu,\bPsi) = \frac{n}{2}\, (2-\delta_{i,j}) \, \big(\E[(X_i-\mu_i)(X_j-\mu_j)] - a_{i,j}\big),
$$
$i,j=1,\ldots,d$.  Since
\begin{align*}
\E[(X_i-\mu_i)(X_j-\mu_j)] &= \E[(X_i-\nu_i+\nu_i-\mu_i)(X_j-\nu_j+\nu_j-\mu_j)] \\
&= \E[(X_i-\nu_i)(X_j-\nu_j) + (\nu_i-\mu_i)(\nu_j-\mu_j) \\
&= \lambda_{i,j} + (\nu_i-\mu_i)(\nu_j-\mu_j)
\end{align*}
then
$$
\frac{\partial}{\partial\psi_{i,j}}\tilde\ell(\bmu,\bPsi) = \frac{n}{2}\, (2-\delta_{i,j}) \, \big(\lambda_{i,j} + (\nu_i-\mu_i)(\nu_j-\mu_j) - a_{i,j}\big).
$$
Next we simultaneously set these derivatives equal to $0$, evaluating the expressions at $(\bmu,\bPsi) = (\widehat\bmu,\widehat\bPsi)$.  Denoting by $\widehat{a}_{i,j}$ the $(i,j)$th entry of $\SS(\widehat\bmu)$, we obtain the system of equations
$$
\widehat\lambda_{i,j} + (\widehat\nu_i - \widehat\mu_i)(\widehat\nu_j - \widehat\mu_j) - \widehat{a}_{i,j} = 0
$$
for all $i,j=1,\ldots,d$.  Denoting $\SS(\widehat\bmu)$ by $\widehat\SS$ and writing these equations in matrix form, we obtain
\begin{equation}
\label{eq_score_Lambda_hat}
\widehat\bLambda + (\widehat\bnu - \widehat\bmu)(\widehat\bnu - \widehat\bmu)' - \widehat\SS = \bzero.
\end{equation}
Since $\widehat\bnu = \overline{\xx}$ then standard algebraic manipulations yield
\begin{align}
\label{eqn:A-hat}
\widehat\SS &= n^{-1}\,\sum_{j=1}^n (\xx_j-\widehat\bmu)(\xx_j-\widehat\bmu)' \nonumber\\
&=n^{-1} \sum_{j=1}^n (\xx_j-\widehat\bnu + \widehat\bnu-\widehat\bmu)(\xx_j-\widehat\bnu + \widehat\bnu-\widehat\bmu)' \nonumber \\
&= n^{-1}\sum_{j=1}^n (\xx_j-\widehat\bnu)(\xx_j-\widehat\bnu)'  +  (\widehat\bnu-\widehat\bmu)(\widehat\bnu-\widehat\bmu)'.
\end{align}
Substituting this result into \eqref{eq_score_Lambda_hat}, we find that
$$
\widehat\bLambda = \widehat\SS - (\widehat\bnu - \widehat\bmu)(\widehat\bnu - \widehat\bmu)' = n^{-1} \sum_{j=1}^n (\xx_j-\overline{\xx})(\xx_j-\overline{\xx})' = \SS(\overline{\xx}).
$$
Finally, by substituting for $\widehat\bLambda$ from \eqref{eq_Lambda_MLE}, we obtain \eqref{eq_score_eqs_Lambda}.
$\qed$

\bigskip

\noindent\textit{Proof of Proposition \ref{prop_necessary_condition}}.
Let $\UU$ be a positive definite $d \times d$ matrix and $\vv \in \mathbb{R}^d$.  By Woodbury's theorem \cite[p.~428, Eq.~(2.25)]{Harville},
$$
(\UU + \vv\vv')^{-1} = \UU^{-1} - \frac{\UU^{-1}\vv\vv'\UU^{-1}}{1+\vv'\UU^{-1}\vv}.
$$
Multiplying this identity on the left by $\vv'$ and on the right by $\vv$, and simplifying the result, we obtain
\begin{align*}
\vv'(\UU + \vv\vv')^{-1}\vv &= \vv'\UU^{-1}\vv - \frac{(\vv'\UU^{-1}\vv)^2}{1+\vv'\UU^{-1}\vv} \\
&= \frac{(1+\vv'\UU^{-1}\vv)\vv'\UU^{-1}\vv - (\vv'\UU^{-1}\vv)^2}{1+\vv'\UU^{-1}\vv}
= \frac{\vv'\UU^{-1}\vv}{1+\vv'\UU^{-1}\vv}.
\end{align*}
This proves that $\vv'(\UU + \vv\vv')^{-1}\vv \in [0,1)$.

Now suppose that the score equations \eqref{eq_score_eqs_nu} and \eqref{eq_score_eqs_Lambda} have a solution.  Setting $\UU = n^{-1} \sum_{j=1}^n (\xx_j-\overline{\xx})(\xx_j-\overline{\xx})'$, and noting that $\UU$ is positive definite, almost surely, for $n > d$, and also setting $\vv = \overline{\xx}-\widehat{\bmu}$, it follows from \eqref{eqn:A-hat} that $\UU + \vv\vv' = \widehat\SS$.  Therefore $(\overline{\xx}-\widehat{\bmu})'\widehat\SS^{-1}(\overline{\xx}-\widehat{\bmu}) \equiv \vv'(\UU + \vv\vv')^{-1}\vv \in [0,1)$.
$\qed$

\bigskip

\noindent\textit{Proof of Proposition \ref{prop_domain_means_upper}}.  
For any $\vv \in \bR^d$, 
\begin{equation}
\label{eq_domain_means_up_1}
\vv'\bLambda\vv = \E \big(\vv'(\XX\XX' - \bnu\bnu')\vv\big) \le \E |\vv'(\XX\XX' - \bnu\bnu')\vv|.
\end{equation}
Since, for any symmetric $d \times d$ matrix $\AA$, there holds $|\vv'\AA\vv| \le \|\vv\|^2 \|\AA\|_F$, 
then 
\begin{align*}
|\vv'(\XX\XX' - \bnu\bnu')\vv| &\le \|\vv\|^2 \|\XX\XX' - \bnu\bnu'\|_F \\
&\le \|\vv\|^2 (\|\XX\XX'\|_F + \|\bnu\bnu'\|_F) 
= \|\vv\|^2 (\|\XX\|^2 + \|\bnu\|^2).
\end{align*}
Applying the latter inequality in \eqref{eq_domain_means_up_1}, we find that 
$$
\vv'\bLambda\vv \le \|\vv\|^2 (\E \|\XX\|^2 + \|\bnu\|^2),
$$
proving that $\bLambda \preccurlyeq (\E (\|\XX\|^2) + \|\bnu\|^2) \II_d$.  Substituting $\bLambda = \E(\XX\XX') - \bnu\bnu'$ and then simplifying, we obtain 
$$
\E(\XX\XX') \preccurlyeq \bnu\bnu' + (\E (\|\XX\|^2) + \|\bnu\|^2) \II_d,
$$
as stated in \eqref{eq_domain_means_upper_1}.  Finally, \eqref{eq_domain_means_upper_2} is obtained by using the fact that $\E(\|\XX\|^2) = \tr\E(\XX\XX')$ and substituting $(\btheta,\bTheta)$ for $(\E(\XX),\E(\XX\XX'))$ in \eqref{eq_domain_means_upper_1}.  
$\qed$

\subsection{Proofs for Section \ref{sec_several_untruncated}}
\label{subsec_sec_several_untruncated}

\noindent\textit{Proof of Proposition \ref{prop_space_single}}.
On writing
$$
\btheta'\tt - \tt' \bTheta \tt = \theta_1 t_1 + \btheta_2' \tt_2 - \Theta_{11} t_1^2 - 2 t_1 \bTheta_{21}' \tt_2 - \tt_2' \bTheta_{22} \tt_2,
$$
it follows that
\begin{multline*}
L_T(\btheta,\bTheta) \\
= \int_0^\infty \exp(\theta_1 t_1 - \Theta_{11} t_1^2) \int_{\bR^{d-1}} \exp((\btheta_2 - 2 t_1 \bTheta_{21})'\tt_2 - \tt_2'\bTheta_{22}\tt_2) \, \dd\tt_2 \dd t_1.
\end{multline*}
For this repeated integral to be finite it follows from Fubini's theorem that the inner integral necessarily converges $t_1$-almost everywhere, a condition that holds if and only if $\bTheta_{22}$ is positive definite.  Applying the multivariate normal integral, i.e.,
\begin{multline*}
\int_{\bR^{d-1}} \exp((\btheta_2-2 t_1 \bTheta_{21})'\tt_2 - \tt_2'\bTheta_{22}\tt_2) \, \dd\tt_2 \\
= \pi^{(d-1)/2} (\det \bTheta_{22})^{-1/2} \exp\big(\tfrac14(\btheta_2-2 t_1 \bTheta_{21})'\bTheta_{22}^{-1}(\btheta_2-2 t_1 \bTheta_{21})\big),
\end{multline*}
we obtain
\begin{align*}
L_T(\btheta,\bTheta) &= \pi^{(d-1)/2} (\det \bTheta_{22})^{-1/2} \\ 
& \ \ \times \int_0^\infty \exp\big(\theta_1 t_1 - \Theta_{11} t_1^2 + \tfrac14(\btheta_2-2 t_1 \bTheta_{21})'\bTheta_{22}^{-1}(\btheta_2-2 t_1 \bTheta_{21})\big) \dd t_1.
\end{align*}

By elementary algebra,
\begin{multline*}
\theta_1 t_1 - \Theta_{11} t_1^2 + \tfrac14(\btheta_2-2 t_1 \bTheta_{21})'\bTheta_{22}^{-1}(\btheta_2-2 t_1 \bTheta_{21}) \\
= -\Theta_{11\cdot 2}t_1^2 + \theta_{1\cdot 2}t_1 + \tfrac14 \btheta_2'\bTheta_{22}^{-1}\btheta_2,
\end{multline*}
so we obtain
\begin{multline}
\label{eq_LT_one_dim_int_sec5}
L_T(\btheta,\bTheta) = \pi^{(d-1)/2} (\det \bTheta_{22})^{-1/2} \exp(\tfrac14 \btheta_2'\bTheta_{22}^{-1}\btheta_2) \\
\times \int_0^\infty \exp(-\Theta_{11\cdot 2}t_1^2 + \theta_{1\cdot 2}t_1) \dd t_1,
\end{multline}
and the latter integral evidently converges if and only if $\Theta_{11\cdot 2} > 0$, or $\Theta_{11\cdot 2} = 0$ and $\theta_{1\cdot 2} < 0$.  Therefore it follows that $L_T(\btheta,\bTheta) < \infty$ if and only if $\bTheta_{22}$ is positive definite and either (i) $\Theta_{11\cdot 2} = 0$ and $\theta_{1\cdot 2} < 0$, or (ii) $\Theta_{11\cdot 2} > 0$.

In the case of (ii), the conditions that $\bTheta_{22}$ is positive definite and $\Theta_{11\cdot 2} > 0$ together are equivalent to the positive definiteness of $\bTheta$.  
To prove this, it suffices to show that $\det(\bTheta) > 0$; however this follows immediately from the Schur complement determinant formula, $\det \bTheta = \Theta_{11 \cdot 2} \cdot \det \bTheta_{22}$.

Therefore $\mcalD = \mcalD_0 \cup \mcalD_1$ where
\begin{equation}
\label{eq_D_0_space}
\mcalD_0 = \{(\btheta,\bTheta) \in \bR^d \times \mcalS^{d \times d}: \theta_{1\cdot 2} < 0, \bTheta_{22} \ \textrm{is positive definite, and } \Theta_{11\cdot 2} = 0\}
\end{equation}
is the canonical parameter space arising from (i) and $\mcalD_1 = \bR^d \times \mcalS_d^{d \times d}$ is the canonical parameter space corresponding to (ii).
$\qed$

\medskip

\noindent\textit{Proof of Theorem \ref{thm_single_comp_non_steep}}.  
Let $\theta_{0,1} < 0$, $\bTheta_{0,22}$ be a $(d-1)\times (d-1)$ positive definite matrix, and define 
\begin{equation}
\label{eq_omega0}
\btheta_0 = \begin{pmatrix} \theta_{0,1} \\ \bzero \end{pmatrix},
\quad
\bTheta_0 = 
\begin{pmatrix} 0 & \bzero \\ \bzero & \bTheta_{0,22} \end{pmatrix}.
\end{equation}
Consider the parameter pair $\bomega_0 = (\btheta_0,\bTheta_0) \in \bR^d \times \mcalS^{d \times d}$; we now show that $\bomega_0 \in \partial\Omega$, the boundary of $\Omega=\hbox{Int}(\mcalD)$, where $\mcalD$ is the canonical parameter space defined in Proposition \ref{prop_space_single}.  

For $n = 1,2,3,\ldots$, define a sequence $\bomega_n = (\btheta_0,\bTheta_n)$, where
$$
\bTheta_n = \begin{pmatrix} n^{-1} & \bzero \\ \bzero & \bTheta_{0,22} \end{pmatrix}.
$$
By \eqref{eq_D_0_space}, $\bomega_n \in \mcalD_0 = {\rm Int}(\mcalD_0)$, where the latter equality holds because $\mcalD_0$ is open.  It is evident that $\bomega_n \to \bomega_0$ as $n \to \infty$, consequently, $\bomega_0 \in \cl({\rm Int}(\mcalD_1)) \subseteq \cl({\rm Int}(\mcalD)) = \cl(\Omega)$.

Next, with the same $\btheta_0$ and $\bTheta_{0,22}$, define a new sequence $\widetilde\bomega_n = (\btheta_0,\widetilde\bTheta_n)$ where 
$$
\widetilde\bTheta_n = \begin{pmatrix} -n^{-1} & \bzero \\ \bzero & \bTheta_{0,22} \end{pmatrix}.
$$
Since $\btheta_0 \neq \bzero$ and $\widetilde\bTheta_n$ is not positive definite then $L_T(\btheta_n,\widetilde\bTheta_n) = \infty$; therefore $\widetilde\bomega_n \in \mcalD^c \subseteq \Omega^c$.  Also, it is evident that $\widetilde\bomega_n \to \bomega_0$ as $n \to \infty$, so it follows that $\bomega_0 \in \cl(\mcalD^c) \subseteq \cl(\Omega^c)$.  Thereby, $\bomega_0 \in \cl(\Omega) \cap \cl(\Omega^c) = \partial\Omega$.

To establish the non-steepness property, recall that the gradient operator is $\nabla_\bomega = (\nabla_\btheta,\nabla_\bTheta)$, having the norm defined in \eqref{eq_nabla_omega}.  We now calculate $\|\nabla_\bomega K(\btheta,\bTheta)\|$ at $\bomega = \bomega_0$.

By the definition of $\theta_{1\cdot 2}$ in \eqref{eq_conditional_thetas}, it follows that $\nabla_{\theta_1} \theta_{1\cdot 2} = 1$.  Therefore by \eqref{eq_LT_one_dim_int_sec5}, and after interchanging the derivative and the integral, it follows that
\begin{align*}
\nabla_{\theta_1} L_T(\btheta,\bTheta)
&= C_d(\bTheta_{22}) \exp(\tfrac14 \btheta_2'\bTheta_{22}^{-1}\btheta_2) \nabla_{\theta_1} \int_0^\infty \exp(-\Theta_{11\cdot 2}t_1^2 + \theta_{1\cdot 2}t_1) \dd t_1 \\
&= C_d(\bTheta_{22}) \exp(\tfrac14 \btheta_2'\bTheta_{22}^{-1}\btheta_2) \int_0^\infty t_1 \exp(-\Theta_{11\cdot 2}t_1^2 + \theta_{1\cdot 2}t_1) \dd t_1.
\end{align*}
Therefore
\begin{align*}
\nabla_{\theta_1} K_T(\btheta,\bTheta) &= \dfrac{\nabla_{\theta_1} L_T(\btheta,\bTheta)}{L_T(\btheta,\bTheta)}
=\frac{\int_0^\infty t_1 \exp(-\Theta_{11\cdot 2}t_1^2 + \theta_{1\cdot 2}t_1) \dd t_1}{\int_0^\infty \exp(-\Theta_{11\cdot 2}t_1^2 + \theta_{1\cdot 2}t_1) \dd t_1},
\end{align*}
and on evaluating this gradient at $(\btheta_0,\bTheta_0)$, we obtain
\begin{equation}
\label{eq_norm1}
\nabla_{\theta_1} K_T(\btheta,\bTheta)\big|_{(\btheta_0,\bTheta_0)} = \frac{\int_0^\infty t_1 \exp(\theta_{0,1}\,t_1) \dd t_1}{\int_0^\infty \exp(\theta_{0,1}\,t_1) \dd t_1}
= \theta_{0,1}^{-1}.
\end{equation}

Next, we again deduce from \eqref{eq_conditional_thetas} that
$$
\nabla_{\btheta_2} (\theta_{1\cdot 2} t_1 + \tfrac14 \btheta_2'\bTheta_{22}^{-1}\btheta_2) = -\bTheta_{22}^{-1}\bTheta_{21}t_1 + \tfrac12 \bTheta_{22}^{-1}\btheta_2,
$$
and then it follows that
\begin{align*}
\nabla_{\btheta_2} L_T(\btheta,\bTheta)
&= C_d(\bTheta_{22}) \nabla_{\btheta_2} \int_0^\infty \exp(-\Theta_{11\cdot 2}t_1^2 + \theta_{1\cdot 2}t_1 + \tfrac14 \btheta_2'\bTheta_{22}^{-1}\btheta_2) \dd t_1 \\
&= C_d(\bTheta_{22}) \exp(\tfrac14 \btheta_2'\bTheta_{22}^{-1}\btheta_2) \\
& \quad \times \int_0^\infty [-\bTheta_{22}^{-1}\bTheta_{21}t_1 + \tfrac12 \bTheta_{22}^{-1}\btheta_2] \exp(-\Theta_{11\cdot 2}t_1^2 + \theta_{1\cdot 2}t_1) \dd t_1.
\end{align*}
Therefore
\begin{align*}
\nabla_{\btheta_2} K_T(\btheta,\bTheta) &= \dfrac{\nabla_{\btheta_2} L_T(\btheta,\bTheta)}{L_T(\btheta,\bTheta)} \\
&= \frac{\int_0^\infty [-\bTheta_{22}^{-1}\bTheta_{21}t_1 + \tfrac12 \bTheta_{22}^{-1}\btheta_2] \exp(-\Theta_{11\cdot 2}t_1^2 + \theta_{1\cdot 2}t_1) \dd t_1}{\int_0^\infty \exp(-\Theta_{11\cdot 2}t_1^2 + \theta_{1\cdot 2}t_1) \dd t_1},
\end{align*}
and, by \eqref{eq_omega0},
\begin{align}
\label{eq_norm2}
\nabla_{\btheta_2} K_T(\btheta,\bTheta)\big|_{(\btheta_0,\bTheta_0)} = \bzero.
\end{align}

Now we calculate the matrix $\nabla_\bTheta K_T(\btheta,\bTheta)$ and evaluate it at $(\btheta_0,\bTheta_0)$.  By definition,
\begin{equation}
\label{eq_Theta_grad_K_sec5}
\nabla_\bTheta K_T(\btheta,\bTheta) = \dfrac{\nabla_\bTheta L_T(\btheta,\bTheta)}{L_T(\btheta,\bTheta)} = - \frac{\int_0^\infty \int_{\bR^{d-1}} \tt\tt' \exp(\btheta'\tt - \tt' \bTheta \tt) \, \dd\tt_2 \dd t_1}{L_T(\btheta,\bTheta)}.
\end{equation}
Define the $d \times d$ symmetric matrix
\begin{equation}
\label{eq_M_matrix_sec5}
\MM = (m_{i,j}) = \int_0^\infty \int_{\bR^{d-1}} \tt\tt' \exp(\btheta'\tt - \tt' \bTheta \tt) \, \dd\tt_2 \dd t_1.
\end{equation}
Then it follows from \eqref{eq_M_matrix_sec5} that
$$
m_{i,j} = \frac{\partial^2}{\partial\theta_i \partial\theta_j} \int_0^\infty \int_{\bR^{d-1}}  \exp(\btheta'\tt - \tt' \bTheta \tt) \, \dd\tt_2 \dd t_1 = \frac{\partial^2}{\partial\theta_i \partial\theta_j} L_T(\btheta,\bTheta).
$$
By \eqref{eq_LT_one_dim_int_sec5},
\begin{align*}
m_{1,1} &= C_d(\bTheta_{22}) \frac{\partial^2}{\partial\theta_1^2} \int_0^\infty \exp(-\Theta_{11\cdot 2}t_1^2 + \theta_{1\cdot 2}t_1 + \tfrac14 \btheta_2'\bTheta_{22}^{-1}\btheta_2) \dd t_1 \\
&= C_d(\bTheta_{22}) \exp(\tfrac14 \btheta_2'\bTheta_{22}^{-1}\btheta_2) \int_0^\infty t_1^2 \exp(-\Theta_{11\cdot 2}t_1^2 + \theta_{1\cdot 2}t_1) \dd t_1.
\end{align*}

For $j=2,\ldots,d$, a similar calculation shows that
$$
m_{j,j} = C_d(\bTheta_{22}) \frac{\partial^2}{\partial\theta_j^2} \int_0^\infty \exp(-\Theta_{11\cdot 2}t_1^2 + \theta_{1\cdot 2}t_1 + \tfrac14 \btheta_2'\bTheta_{22}^{-1}\btheta_2) \dd t_1.
$$
Noting that $\theta_{1\cdot 2} = \theta_1 - \bTheta_{12}\bTheta_{22}^{-1}\btheta_2 \equiv \theta_1 - \btheta_2'\bTheta_{22}^{-1}\bTheta_{21}$, we have
$$
\nabla_{\btheta_2} (-\Theta_{11\cdot 2}t_1^2 + \theta_{1\cdot 2}t_1 + \tfrac14 \btheta_2'\bTheta_{22}^{-1}\btheta_2) = -\bTheta_{22}^{-1}\bTheta_{21}t_1 + \tfrac12 \bTheta_{22}^{-1}\btheta_2,
$$
equivalently, for $j=2,\ldots,d$,
\begin{align*}
\frac{\partial}{\partial\theta_j} (-\Theta_{11\cdot 2}t_1^2 + \theta_{1\cdot 2}t_1 + \tfrac14 \btheta_2'\bTheta_{22}^{-1}\btheta_2)
&= (-\bTheta_{22}^{-1}\bTheta_{21}t_1 + \tfrac12 \bTheta_{22}^{-1}\btheta_2)_j \\
&\equiv -(\bTheta_{22}^{-1}\bTheta_{21})_j t_1 + \tfrac12 \sum_{k=2}^d (\bTheta_{22}^{-1})_{j,k} \theta_k.
\end{align*}
Therefore, for $j=2,\ldots,d$,
\begin{align*}
\frac{\partial^2}{\partial\theta_j^2} (-\Theta_{11\cdot 2}t_1^2 + \theta_{1\cdot 2}t_1 + \tfrac14 \btheta_2'\bTheta_{22}^{-1}\btheta_2) = \tfrac12 (\bTheta_{22}^{-1})_{j,j}.
\end{align*}

Next, it follows by direct calculation that
\begin{align*}
\frac{\partial^2}{\partial\theta_j^2} \exp(-&\Theta_{11\cdot 2}t_1^2 + \theta_{1\cdot 2}t_1 + \tfrac14 \btheta_2'\bTheta_{22}^{-1}\btheta_2) \\
&= \big[\big((-\bTheta_{22}^{-1}\bTheta_{21}t_1 + \tfrac12 \bTheta_{22}^{-1}\btheta_2)_j\big)^2 + \tfrac12 (\bTheta_{22}^{-1})_{j,j}\big] \\
&\qquad\qquad\quad \times \exp(-\Theta_{11\cdot 2}t_1^2 + \theta_{1\cdot 2}t_1 + \tfrac14 \btheta_2'\bTheta_{22}^{-1}\btheta_2),
\end{align*}
hence
\begin{align*}
\sum_{j=2}^d \frac{\partial^2}{\partial\theta_j^2} & \exp(-\Theta_{11\cdot 2}t_1^2 + \theta_{1\cdot 2}t_1 + \tfrac14 \btheta_2'\bTheta_{22}^{-1}\btheta_2) \\
&= \bigg(\sum_{j=2}^d \big[\big((-\bTheta_{22}^{-1}\bTheta_{21}t_1 + \tfrac12 \bTheta_{22}^{-1}\btheta_2)_j\big)^2 + \tfrac12 (\bTheta_{22}^{-1})_{j,j}\big]\bigg) \\
&\qquad\qquad\qquad \times \exp(-\Theta_{11\cdot 2}t_1^2 + \theta_{1\cdot 2}t_1 + \tfrac14 \btheta_2'\bTheta_{22}^{-1}\btheta_2).
\end{align*}
Since
\begin{align*}
\sum_{j=2}^d & \big[\big((-\bTheta_{22}^{-1}\bTheta_{21}t_1 + \tfrac12 \bTheta_{22}^{-1}\btheta_2)_j\big)^2 + \tfrac12 (\bTheta_{22}^{-1})_{j,j}\big] \\
&= (-\bTheta_{22}^{-1}\bTheta_{21}t_1 + \tfrac12 \bTheta_{22}^{-1}\btheta_2)'(-\bTheta_{22}^{-1}\bTheta_{21}t_1 + \tfrac12 \bTheta_{22}^{-1}\btheta_2) + \tfrac12 \tr(\bTheta_{22}^{-1}) \\
&= \bTheta_{12}\bTheta_{22}^{-2}\bTheta_{21}t_1^2 + \tfrac14 \btheta_2'\bTheta_{22}^{-2}\btheta_2 - \bTheta_{12}\bTheta_{22}^{-2}\btheta_2 t_1 + \tfrac12 \tr(\bTheta_{22}^{-1}),
\end{align*}
then
\begin{align*}
\sum_{j=2}^d m_{j,j} &= C_d(\bTheta_{22}) \sum_{j=2}^d \frac{\partial^2}{\partial\theta_j^2} \int_0^\infty \exp\big(-\Theta_{11\cdot 2}t_1^2 + \theta_{1\cdot 2}t_1 + \tfrac14 \btheta_2'\bTheta_{22}^{-1}\btheta_2\big) \dd t_1 \\
&= C_d(\bTheta_{22}) \\
&\quad \times \int_0^\infty \big[\bTheta_{12}\bTheta_{22}^{-2}\bTheta_{21}t_1^2 + \tfrac14 \btheta_2'\bTheta_{22}^{-2}\btheta_2 - \bTheta_{12}\bTheta_{22}^{-2}\btheta_2 t_1 + \tfrac12 \tr(\bTheta_{22}^{-1})\big] \\
&\qquad\qquad\qquad\qquad\qquad\quad \times \exp(-\Theta_{11\cdot 2}t_1^2 + \theta_{1\cdot 2}t_1 + \tfrac14 \btheta_2'\bTheta_{22}^{-1}\btheta_2) \dd t_1.
\end{align*}
Collecting these results together, we have obtained
\begin{align*}
\tr&(\MM) = \sum_{j=1}^d m_{j,j} \\
&= C_d(\bTheta_{22}) \exp(\tfrac14 \btheta_2'\bTheta_{22}^{-1}\btheta_2) \\
&\quad \times \int_0^\infty \big[(1+\bTheta_{12}\bTheta_{22}^{-2}\bTheta_{21})t_1^2 + \tfrac14 \btheta_2'\bTheta_{22}^{-2}\btheta_2 - \bTheta_{12}\bTheta_{22}^{-2}\btheta_2 t_1 + \tfrac12 \tr(\bTheta_{22}^{-1})\big] \\
&\qquad\qquad\qquad\qquad\qquad\qquad\qquad\qquad\qquad \times \exp(-\Theta_{11\cdot 2}t_1^2 + \theta_{1\cdot 2}t_1) \dd t_1,
\end{align*}
therefore
\begin{align}
\label{eq_tr_M}
\tr(\MM)\big|_{(\btheta_0,\bTheta_0)} &= C_d(\bTheta_{22}) \int_0^\infty \big[t_1^2 + \tfrac12 \tr(\bTheta_{22}^{-1})\big] \exp(\theta_{0,1}\,t_1) \dd t_1 \notag \\
&= -\big(2\theta_{0,1}^{-3} + \tfrac12 \tr(\bTheta_{22}^{-1}) \theta_{0,1}\big) C_d(\bTheta_{22}).
\end{align}

Denote by $\eta_1(\MM),\ldots,\eta_d(\MM)$ the eigenvalues of $\MM$.  Also, let $\|\MM\|_1 := \sum_{j=1}^d |\eta_j(\MM)|$ denote the $\ell_1$-norm of $\MM$.  The Frobenius norm, $\|\MM\|_F$, of $\MM$ satisfies
$$
\|\MM\|_F = [\tr(\MM^2)]^{1/2} = \bigg(\sum_{j=1}^d [\eta_j(\MM)]^2\bigg)^{1/2} \le \sum_{j=1}^d |\eta_j(\MM)| = \|\MM\|_1.
$$

By \eqref{eq_M_matrix_sec5}, $\MM$ also is positive semidefinite, hence $\eta_j(\MM) \ge 0$, $j=1,\ldots,d$.  Therefore $\|\MM\|_1 = \tr(\MM)$, and it follows from \eqref{eq_Theta_grad_K_sec5} that
$$
\|\nabla_\bTheta K_T(\btheta,\bTheta)\|_F = \frac{\|\MM\|_F}{L_T(\btheta,\bTheta)} \le \frac{\|\MM\|_1}{L_T(\btheta,\bTheta)} = \frac{\tr(\MM)}{L_T(\btheta,\bTheta)}.
$$
By \eqref{eq_tr_M},
\begin{multline}
\|\nabla_\bTheta K_T(\btheta,\bTheta)\|_F \Big|_{(\btheta_0,\bTheta_0)}  \\
\le \frac{\tr(\MM)\big|_{(\btheta_0,\bTheta_0)}}{L_T(\btheta_0,\bTheta_0)}= -\frac{\big(2\theta_{0,1}^{-3} + \tfrac12 \tr(\bTheta_{22}^{-1}) \theta_{0,1}\big) C_d(\bTheta_{22})}{L_T(\btheta_0,\bTheta_0)}.
\label{eq_norm3}
\end{multline}
Combining \eqref{eq_norm1}, \eqref{eq_norm2}, and \eqref{eq_norm3}, we obtain
\begin{align*}
    \|\nabla_\bomega K(\btheta,\bTheta)\|^2\big|_{(\btheta_0,\bTheta_0)}  &= \|\nabla_\btheta K(\btheta,\bTheta)\|^2\big|_{(\btheta_0,\bTheta_0)}  + \|\nabla_\bTheta K(\btheta,\bTheta)\|_F^2\big|_{(\btheta_0,\bTheta_0)}\\
    &\leq \theta_{0,1}^{-2} + \bigg[\frac{\big(2\theta_{0,1}^{-3} + \tfrac12 \tr(\bTheta_{22}^{-1}) \theta_{0,1}\big) C_d(\bTheta_{22})}{L_T(\btheta_0,\bTheta_0)}\bigg]^2,
\end{align*}
which is finite.
Therefore, as $n \to \infty$, 
$$
\|\nabla_\bomega K(\btheta,\bTheta)\|^2 \Big|_{(\btheta_n,\bTheta_n)} \to \|\nabla_\bomega K(\btheta,\bTheta)\|^2\Big|_{(\btheta_0,\bTheta_0)} < \infty,
$$
proving that the exponential family corresponding to \eqref{eq_LT_1_truncated_sec5} is not steep.
$\qed$

\bigskip
\bigskip

\noindent
\textbf{Acknowledgments}.  We are very grateful to the reviewers for their supremely careful reading of the manuscript and for providing us with incisive and detailed comments that led to numerous improvements.  

\medskip

\noindent
\textbf{Declarations}.  {The research of Michael Levine is partially supported by NSF-DMS grant \# 2311103. The authors have no relevant financial or non-financial interests to disclose. The authors have no conflicts of interest to declare.}


\bigskip

\begin{thebibliography}{999}
\bibliographystyle{ims}

{\parskip=3pt

\bibitem{Amemiya}
Amemiya, T. (1974).  Multivariate regression and simultaneous equations models when the dependent variables are truncated normal.  \textit{Econometrica}, {\bf 42}, 999-1012.

\bibitem{Anderson}
Anderson, T. W. (2003).  \textsl{An Introduction to Multivariate Statistical Analysis}, third edition.  Wiley, New York, NY.

\bibitem{Bar_Lev}
Bar-Lev, S. K., and Enis, P. (1986).  Reproducibility and natural exponential families with power variance functions.  \textit{Ann. Statist.}, {\bf 14}, 1507--1522.

\bibitem{Bar_Lev_Ridder}
Bar-Lev, S. K., and Ridder, A. (2021).  New exponential dispersion models for count data: the ABM and LM classes.  \textit{ESAIM: Probab. Statist.}, {\bf 25}, 31--52.

\bibitem{Barndorff}
Barndorff-Nielsen, O. E. (1978).  \textsl{Information and Exponential Families in Statistical Theory}, (second printing 2014). Wiley, Chichester, U.K.

\bibitem{Bedbur_etal}
Bedbur, S., Kamps, U., and Imm, A. (2023).  On the existence of maximum likelihood estimates for the parameters of the Conway-Maxwell-Poisson distribution.  \textit{ALEA, Lat. Am. J. Probab. Math. Stat.}, {\bf 20}, 561--575.

\bibitem{Brown}
Brown, L. D. (1986).  \textsl{Fundamentals of Statistical Exponential Families: With Applications in Statistical Decision Theory}, Institute of Mathematical Statistics, Hayward, CA.

\bibitem{Burkill}
Burkill, J. C., and Burkill, H. (2002).  \textsl{A Second Course in Mathematical Analysis}. Cambridge University Press, New York.

\bibitem{Cohen}
Cohen, A. C. (2016).  \textsl{Truncated and Censored Samples: Theory and Applications}, CRC Press, Boca Raton, FL.

\bibitem{delCastillo}
del Castillo, J. (1994).  The singly truncated normal distribution: A non-steep exponential family.  \textit{Ann. Inst. Statist. Math.}, \textbf{46}, 57--66.

\bibitem{Gaal}
Gaal, S. A. (2009).  \textsl{Point Set Topology}. Dover Publications, Mineola, NY.

\bibitem{Harville}
Harville, D. A. (2008).  \textsl{Matrix Algebra From a Statistician's Perspective}, Springer, New York, NY.

\bibitem{Hegde}
Hegde, L. M., and Dahiya, R. C. (1989).  Estimation of the parameters in a truncated normal distribution.  \textit{Commun. Statist. - Theory Methods}, \textbf{18}, 4177--4195.

\bibitem{Heckman}
Heckman, J. J. (1976).  The common structure of statistical models of truncation, sample selection and limited dependent variables and a simple estimator for such models. \textit{Ann. Econ. Social Measurement}, \textbf{5}, 475--492.

\bibitem{Hong}
Hong, H., and Shum, M. (2002).  Econometric models of asymmetric ascending auctions. \textit{J. Econometrics}, \textbf{112}, 327--358.

\bibitem{Horn}
Horn, R. A., and Johnson, C. R. (1990).  \textsl{Matrix Analysis}. Cambridge University Press, New York, NY.

\bibitem{Jin}
Jin, B. S., Han, J. J., Ding, S., and Miao, B. Q. (2018).  EM algorithm of the truncated multinormal distribution with linear restriction on the variables. \textit{Acta Math. Appl. Sinica}, \textbf{34,} 155--162.

\bibitem{Kan}
Kan, R., and Robotti, C. (2017).  On moments of folded and truncated multivariate normal distributions. \textit{J. Comput. Graphical Statist.}, \textbf{26}, 930--934.

\bibitem{Lachos}
Lachos, V. H., Matos, L. A., Castro, L. M., and Chen, M.-H. (2019).  Flexible longitudinal linear mixed models for multiple censored responses data. \textit{Statist. Med.}, \textbf{38}, 1074--1102.

\bibitem{Leppard}
Leppard, P., and Tallis, G. M. (1989).  Algorithm AS 249: Evaluation of the mean and covariance of the truncated multinormal distribution.  \textit{J. Roy. Statist. Soc.}, Ser. C, \textbf{38}, 543--553.

\bibitem{Levine}
Levine, M., Richards, D., and Su, J. (2020).  Independence properties of the truncated multivariate elliptical distributions. \textit{Statist. Probab. Lett.}, \textbf{161}, 108729.

\bibitem{Lin}
Lin, T. I., and Wang, W. L. (2017).  Multivariate-$t$ nonlinear mixed models with application to censored multi-outcome AIDS studies. \textit{Biostatistics}, \textbf{18}, 666--681.

\bibitem{Malley}
Malley, J. D. (1983).  Statistical and algebraic independence.  \textit{Ann. Statist.}, \textbf{11}, 341--345.

\bibitem{Manjunath}
Wilhelm, S., and Manjunath, B. G. (2010).  {tmvtnorm}: A package for the truncated multivariate normal distribution.  \textit{The R Journal}, \textbf{2}, 25--29.

\bibitem{Rosenbaum}
Rosenbaum, S. (1961).  Moments of a truncated bivariate normal distribution. \textit{J. Roy. Statist. Soc.}, Ser. B, \textbf{23}, 405--408.

\bibitem{Sundberg}
Sundberg, R. (2019).  \textsl{Statistical Modelling by Exponential Families}.  Cambridge University Press, New York.

\bibitem{Tallis}
Tallis, G. M. (1961).  The moment-generating function of the truncated multi-normal distribution. \textit{J. Roy. Statist. Soc.}, Ser. B, \textbf{23}, 223--229.

\bibitem{Yu}
Yu, J.-W., and Tian, G.-L. (2011). Efficient algorithms for generating truncated multivariate normal distributions. \textit{Acta Math. Appl. Sinica}, \textbf{27}, 601--612.

}\end{thebibliography}

\end{document}